\documentclass[11pt]{article}
\usepackage[leqno]{amsmath}   
\usepackage{amsthm}   
\usepackage{amssymb}  
\usepackage{amscd}   
\usepackage{cdextra} 
\begin{document}
\newcommand{\proLog}{\operatorname{{\cal L}og}}
\newcommand{\Log}{\operatorname{{\cal L}og}}
\newcommand{\Logeom}{\operatorname{{\cal L}og}^{\geom}}
\newcommand{\Rhgeom}{\Rh^{\geom}}
\newcommand{\Pic}{\operatorname{P}}
\newcommand{\Alb}{\operatorname{Alb}}
\newcommand{\can}{\operatorname{can}}
\newcommand{\bad}{\operatorname{bad}}
\newcommand{\geom}{\operatorname{g}}
\newcommand{\Ulisse}{\operatorname{U-lisse}}
\newcommand{\CH}{{C_H}}
\newcommand{\CK}{{C_K}}
\newcommand{\J}{\operatorname{P}^0}
\newcommand{\ab}{\operatorname{ab}}
\newcommand{\op}{\operatorname{op}}
\newcommand{\pro}{\operatorname{pro}}
\newcommand{\af}{\frak{a}}
\newcommand{\mf}{\frak{m}}
\newcommand{\bff}{\frak{b}}
\newcommand{\qf}{\frak{q}}
\newcommand{\Uf}{\frak{U}}
\newcommand{\If}{\frak{I}}
\newcommand{\Jf}{\frak{J}}
\newcommand{\Tr}{\operatorname{Tr}}
\newcommand{\CHbar}{\overline{\CH}}
\newcommand{\CHtilde}{\widetilde{C_H}}
\newcommand{\CHpol}{{C_H^{\pol}}}
\newcommand{\Cpol}{{C^{\pol}}}
\newcommand{\CHprime}{{C_{H'}}}
\newcommand{\CHpolprime}{{C_{H'}^{\pol}}}
\newcommand{\shH}{{\cal H}}
\newcommand{\shK}{{\cal K}}
\newcommand{\potshH}[1]{\Sympot^{#1}\shH}
\newcommand{\envelopshHs}{\widehat{\frak{U}}_{\shH_s}}
\newcommand{\Sh}{\operatorname{Sh}}

\newcommand{\car}{\operatorname{char}}
\newcommand{\pol}{\operatorname{{\cal P}ol}}
\newcommand{\polgeom}{\operatorname{{\cal P}ol}^{\geom}}
\newcommand{\rmpol}{\operatorname{Pol}}
\newcommand{\rmlog}{\operatorname{Log}}
\newcommand{\stand}{\operatorname{stand}}
\newcommand{\class}{\operatorname{class}}
\newcommand{\sing}{\operatorname{sing}}
\newcommand{\cris}{\operatorname{cris}}
\newcommand{\ord}{\operatorname{ord}}
\newcommand{\constr}{\operatorname{constr}}
\newcommand{\End}{\operatorname{End}}
\newcommand{\Fr}{\operatorname{Fr}}
\newcommand{\elliptic}{\operatorname{ell}}
\newcommand{\Qbar}{\overline{\Q}}
\newcommand{\Kbar}{\overline{K}}
\newcommand{\Ehbar}{\overline{\Eh}}
\newcommand{\Chbar}{\overline{\Ch}}
\newcommand{\Kinfty}{K_{\infty}}
\newcommand{\Ohinfty}{\Oh_{\infty}}
\newcommand{\Ahinfty}{\Ah_{\infty}}
\newcommand{\Ehbarinfty}{\Ehbar_{\infty}}
\newcommand{\Chbarinfty}{\Chbar_{\infty}}
\newcommand{\Uhinfty}{\Uh_{\infty}}
\newcommand{\Xhinfty}{\Xh_{\infty}}
\newcommand{\Yhinfty}{\Yh_{\infty}}
\newcommand{\Minfty}{M_{\infty}}

\newcommand{\psibar}{\overline{\psi}}
\newcommand{\shP}{\pol}
\newcommand{\vardelta}{\delta}
\newcommand{\ff}{\frak{f}}
\newcommand{\gf}{\frak{g}}
\newcommand{\pf}{\frak{p}}
\newcommand{\Gen}{\operatorname{{\cal G}en}}
\newcommand{\Gext}{\operatorname{{\cal P}ol}}
\newcommand{\Ab}{{\cal A}b}
\newcommand{\Ph}{{\cal P}}
\newcommand{\Ah}{{\cal A}}
\newcommand{\Bh}{{\cal B}}
\newcommand{\Xh}{{\cal X}}
\newcommand{\Uh}{{\cal U}}
\newcommand{\Ch}{{\cal C}}
\newcommand{\Rh}{{\cal R}}
\newcommand{\Oh}{{\cal O}}
\newcommand{\Dh}{{\cal D}}
\newcommand{\Mh}{{\cal M}}
\newcommand{\Eh}{{\cal E}}
\newcommand{\Fh}{{\cal F}}
\newcommand{\Gh}{{\cal G}}
\newcommand{\Jh}{{\cal J}}
\newcommand{\Kh}{{\cal K}}
\newcommand{\Th}{{\cal T}}
\newcommand{\Ih}{{\cal I}}
\newcommand{\Qh}{{\cal Q}}
\newcommand{\Yh}{{\cal Y}}
\newcommand{\Ihgeom}{{\cal I}^{\geom}}
\newcommand{\1}{{\bf 1}}
\newcommand{\et}{}
\newcommand{\incl}{{\operatorname{incl}}}
\newcommand{\ch}{{\operatorname{ch}}}
\newcommand{\Jac}{{\operatorname{Jac}}}
\renewcommand{\P}{{\Bbb{P}}}
\newcommand{\Ptilde}{\widetilde{{\cal P}}}
\newcommand{\Ltilde}{\widetilde{{\cal L}}}
\newcommand{\shF}{{\cal F}}
\newcommand{\shG}{{\cal G}}
\newcommand{\shE}{{\cal E}}
\newcommand{\Hcont}{H_{\operatorname{et}}}
\newcommand{\Het}{H_{\operatorname{et}}}  
\newcommand{\hm}{H_{\cal M}}
\newcommand{\Lh}{{\cal L}}
\newcommand{\Sdach}{\hat{S}}
\newcommand{\m}{M}
\newcommand{\mquer}{\overline{\m}}
\newcommand{\cusp}{\operatorname{Cusp}}
\newcommand{\base}{B}
\newcommand{\Dir}{\operatorname{Dir}}
\newcommand{\Div}{\operatorname{Div}}
\newcommand{\triv}{\operatorname{triv}}

\newcommand{\Ebar}{\overline{E}}
\newcommand{\Cbar}{\overline{C}}
\newcommand{\Xbar}{\overline{X}}
\newcommand{\fbar}{\overline{f}}
\newcommand{\sbar}{\overline{s}}
\newcommand{\xbar}{\overline{x}}
\renewcommand{\tilde}{\widetilde}
\newcommand{\pibar}{\overline{\pi}}
\newcommand{\jint}{j_{!*}}
\newcommand{\shEis}{\operatorname{{\cal E}}}
\newcommand{\kerres}{\Q[E[N]\ohne 0]^0[\res_\infty]}
\newcommand{\tr}{\operatorname{tr}} 
\newcommand{\klammerk}{{(k)}}
\newcommand{\Symb}{\operatorname{Symb}}
\newcommand{\Norm}{\operatorname{N}}
\newcommand{\ad}{\operatorname{ad}}
\newcommand{\fra}{\frak{a}}
\newcommand{\pr}{\operatorname{pr}}
\newcommand{\cl}{\operatorname{cl}}
\newcommand{\tra}{\operatorname{tr}_{[a]}}
\newcommand{\hmot}{H_{{\cal M}}}
\newcommand{\Hmot}{H^{{\cal M}}}
\renewcommand{\Bbb}{\mathbb}  
\newcommand{\Lbb}{{\Bbb{L}}}
\newcommand{\N}{{\Bbb{N}}} 
\newcommand{\Q}{{\Bbb{Q}}}  
\newcommand{\R}{{\Bbb{R}}}  
\newcommand{\C}{{\Bbb{C}}}  
\newcommand{\Z}{{\Bbb{Z}}}  
\newcommand{\Pe}{{\Bbb{P}}} 
\newcommand{\A}{{\Bbb{A}}}  
\newcommand{\G}{{\Bbb{G}}}  
\newcommand{\D}{{\Bbb{D}}}
\newcommand{\F}{{\Bbb{F}}}  
\newcommand{\Gm}{\G_{m}}      
\newcommand{\U}{\Bbb{U}}    
\renewcommand{\epsilon}{\varepsilon}
\renewcommand{\rho}{\varrho}
\renewcommand{\bar}{\overline}

\newcommand{\tensor}{\otimes}   
\newcommand{\isom}{\cong}       
\newcommand{\ohne}{\smallsetminus}
\newcommand{\argdot}{{\;\cdot\;}}
\newcommand{\Hom}{\operatorname{Hom}}
\newcommand{\Ext}{\operatorname{Ext}}
\newcommand{\Kern}{\operatorname{Ker}}
\newcommand{\bild}{\operatorname{Im}}
\newcommand{\prolim}{\varprojlim}
\newcommand{\indlim}{\varinjlim}
\newcommand{\Ind}{\operatorname{Ind}}


\newcommand{\bew}{\begin{proof}}
\newcommand{\bewende}{\end{proof}}
\newtheorem{lemma}{Lemma}[subsection]
\newtheorem{prop}[lemma]{Proposition}
\newtheorem{thm}[lemma]{Theorem}
\newtheorem{defn}[lemma]{Definition}
\newtheorem{Conditions}[lemma]{Conditions}
\newtheorem{assumption}[lemma]{Assumption}
\newtheorem{notation}[lemma]{Notation}
\newtheorem{cor}[lemma]{Corollary}

\newcommand{\bem}{\noindent{\bf Remark:\ }}
\newcommand{\rem}{\noindent{\bf Remark:\ }}
\newtheorem{conj}[lemma]{Conjecture}


\newcommand{\Spec}{\operatorname{Spec}} 
\newcommand{\punkt}{{\bf .}}
\newcommand{\cone}{\operatorname{Cone}}
\newcommand{\Gr}{\operatorname{Gr}}
\newcommand{\gr}{\Gr}
\newcommand{\coh}{\operatorname{coh}}
\newcommand{\Gl}{\operatorname{Gl}}
\newcommand{\fH}{\frak{H}}
\newcommand{\sym}{\frak{S}}
\newcommand{\Sympot}{\operatorname{Sym}}
\newcommand{\Sym}{\operatorname{Sym}}
\newcommand{\envelop}{\widehat{\frak{U}}_{\shH}}
\newcommand{\envelopQl}{\widehat{\frak{U}}_{\Q_l}}
\newcommand{\envelopR}{\widehat{\frak{U}}_{\R}}
\newcommand{\mult}{\operatorname{mult}}

\newcommand{\res}{\operatorname{res}}

\newcommand{\GL}{\operatorname{GL}}
\newcommand{\SL}{\operatorname{SL}}
\newcommand{\Gal}{\operatorname{Gal}}

\newcommand{\Isom}{\operatorname{Isom}}
\newcommand{\Aut}{\operatorname{Aut}}

\newcommand{\tei}{\, | \,}

\newcommand{\Res}{\operatorname{Res}}
\newcommand{\Eis}{\operatorname{Eis}}

\newcommand{\hullet}{{\scriptstyle\bullet}}

\newcommand{\id}{\operatorname{id}}
\newcommand{\sgn}{\operatorname{sgn}}
\newcommand{\verk}{{\scriptstyle\circ}}

\newcommand{\kerr}{\operatorname{ker}}

\setcounter{section}{0}
\begin{center}{\bf \Large The Tamagawa number conjecture for CM elliptic curves
}\\\vspace{0.5cm}
by Guido Kings\\\vspace{0.5cm}
{\large Preprint  March 17, 2000}
\end{center}
\section*{Introduction}
In this paper   we  prove the Tamagawa number conjecture of
Bloch-Kato for CM elliptic curves using a new  explicit description of 
the specialization  of the  elliptic polylogarithm. 

The   Tamagawa number
conjecture describes the special values of the L-function
of a variety in terms of the regulator maps of the K-theory of the variety
into Deligne and \'etale cohomology (see \ref{blochkatoconj}
for the exact formulation). There are only two cases proven so far in the 
non critical situation, both due to Bloch and Kato \cite{Bl-Ka}:
The first is the Riemann zeta function (i.e. the case of $\Q$) and the second
the L-value  at $2$ of a CM elliptic curve defined over $\Q$ for regular primes. 
In the last case Bloch and Kato use an ad hoc method to describe the $p$-adic
regulator of the K-theory. This does not extend to higher K-groups.

The regulator map to Deligne cohomology was computed by Deninger
\cite{Den1} with the help of the Eisenstein symbol. 
Here,  the regulator 
can be described in terms of real analytic Eisenstein series
(whence the name) and leads to a proof of the Beilinson
conjecture for CM elliptic curves. For the Tamagawa number conjecture
one needs an understanding of the $p$-adic regulator on the 
subspace of K-theory defined by the Eisenstein symbol. In an
earlier paper \cite{Hu-Ki} we established with A. Huber the relation of the
$p$-adic regulator of the
Eisenstein symbol with the specialization of the $p$-adic elliptic
polylogarithm sheaf. The problem remains, to compute these
specializations.

The elliptic polylogarithm is one of the most powerful tools
 in the study of special values of L-functions. All
known cases of the Beilinson conjecture are proved or can be proved
with specializations of the elliptic polylog. A universal property 
characterizes the polylog, which simplifies the explicit
computations.   So far, only the absolute Hodge realization
of the elliptic polylog was understood, due to the extensive work of 
Beilinson and Levin \cite{Be-Le}.
Missing was a theory of the $p$-adic realization, which
 give manageable \'etale cohomology classes. Such an explicit realization
is known in the cyclotomic case.
The approach there, mainly due to Deligne \cite{Del1},
uses torsors over
$\Gm\ohne 1$ which are ramified in $1$. This 
is not transferable to the elliptic case, because
such torsors over elliptic curves
do not exist (there are no Galois coverings ramified in 
exactly one point due to compactness). In our approach  we allow instead
ramification at torsion points, constructing in fact torsors
over $E\ohne E[l^n]$. But this is not the only change of point of view compared
to the cyclotomic case. The question is also, what is the group whose torsors
we have to consider. The right choice is the group of 
torsion points of the torus with character group the augmentation ideal 
of the group ring of $E[l^n]$.

It turns out that the elliptic polylogarithm is an inverse 
limit of $p^r$-torsion points of a certain one-motive, which
is essentially the generalized Jacobian 
defined by the divisor of all $p^r$-torsion points on the
elliptic curve $E$. The cohomology classes of the elliptic 
polylogarithm sheaf can then be described by classes of sections
of certain line bundles. These sections are elliptic units and
going carefully through the construction one finds an 
analog of the elliptic Soul\'e elements of \cite{So2}.

Now enters Iwasawa theory: By an idea of Kato, 
going back in part to earlier work of Soul\'e \cite{So2},
the \'etale cohomology groups can
be described in terms of Iwasawa modules and Rubin's ``main conjecture''
\cite{Ru1} allows to give a bound on the kernel and the cokernel of Soul\'e's
map from elliptic units to the \'etale cohomology. 
On the way we also need some of the tools developed by Rubin to prove
the main conjecture and the whole theory is the second decisive input
into the proof of the Tamagawa number conjecture.\\

Let us finally give a rough sketch of the contents of this paper.
More overviews can be found at the beginning of each section.
In the first section we recall
the statement of the Tamagawa number conjecture and formulate our main result \ref{mainthm}.
After this we recall Deninger's construction of the 
K-theory elements leading to the Beilinson conjecture
for CM elliptic curves. Here we also reduce to the computation
of the specialization of the elliptic polylogarithm sheaf by using our
earlier work \cite{Hu-Ki}. 

The second section reviews the ``main conjecture'' and relates
the Iwasawa modules to \'etale cohomology. 

The next two sections are independent of the rest of the paper. The third section
introduces the elliptic polylog and its specializations. The approach 
follows the important paper \cite{Be-Le} but puts the emphasis on 
different aspects, which are important for our geometric construction
of the elliptic polylog. 

The technical heart of the paper is section four. Here we formulate the
polylog as an inverse limit of torsion points of one motives. The cohomology
classes of what we call ``geometric polylog'' are then computed as the
classes of sections of certain line bundles.

The last section puts the various results together and gives the proof of the
main theorem \ref{mainthm}.\\

It is a pleasure to thank Annette Huber for her constant encouragement
 and for 
the many discussions about the contents of this paper and related
results. Pierre Colmez read a first version of this paper and pointed out 
some inaccuracies and gave valuable hints for improvement. I like to thank
him warmly. Thanks go also
to Christopher Deninger for insisting strongly that my computation 
of the $l$-adic elliptic polylog should be turned immediately into 
a proof of the Tamagawa number conjecture. 
Also I like to thank P. Schneider for clarifying some point in
 Iwasawa theory. 

\section{The Bloch-Kato 
conjecture for CM-elliptic curves}

This part of the paper contains our first main result, the Bloch-Kato
conjecture for CM elliptic curves. For the precise formulation of our result
we refer to \ref{mainthm}. This part is organized as follows: First we 
review the Bloch-Kato conjecture in the case of interest to us. 
Here we also formulate the main theorem. Then we recall the construction 
due to Deninger of elements in the K-theory of CM elliptic curves.
These elements satisfy the Beilinson conjecture for these curves as
was shown by Deninger. This is the starting point of our investigations
of the Bloch-Kato conjecture. In the formulation of Kato, the conjecture
is concerned with $\Z_p$-lattices of $\Q_p$-vector spaces namely
\'etale cohomology groups. Using an idea of Soul\'e (the elliptic
Soul\'e elements) and the machinery
of the Iwasawa ``main conjecture'' developed by Rubin, we can describe
the \'etale cohomology groups, or better
the complexes computing them, fairly well. The Iwasawa ``main conjecture''
will be reviewed in  section \ref{secmainconj} and the description of 
the \'etale cohomology groups is in section \ref{secreduction}.
Finally, the last section contains the comparison between
the elliptic Soul\'e elements and Deninger's K-theory elements.
Here we need the results of part \ref{seccomp} and
in particular theorem \ref{polastheta}.

\subsection{The Tamagawa number conjecture of Bloch-Kato and the
main theorem}\label{blochkatoconj}
This section  recalls the Bloch-Kato conjecture \cite{Bl-Ka} about special 
values of L-functions in the formulation of Kato \cite{Ka1} and
\cite{Ka2}. We review
this only for certain weights,
which suffices for our purpose.  Then we formulate our main result.

\subsubsection{The Tamagawa number conjecture in the formulation of 
Kato}\label{tamconj}

Let $X/K$ be a smooth proper variety over a number field
$K$ with ring of integers $\Oh_K$.
Fix integers $m\ge 0$ and $r$ such that $m-2r\le -3$ and $r>\inf(m,\dim(X))$.
Let $p$ be a prime number not equal to $2$.   Let $S$ be a
set of finite primes of $K$ containing
the primes  lying over $p$ and the ones where $X$ has bad reduction.
Let $\Oh_S$ be $\Oh_K[\frac{1}{S}]$ the ring, where the primes
in $S$ are inverted.
Define $\Gal(\Kbar/K)$-modules
\begin{align*}
V_p&:= H^m_{\et}(X\times_{K}\Kbar, \Q_p(r))\\
T_p&:= H^m_{\et}(X\times_{K}\Kbar, \Z_p(r))
\end{align*}
Let $j:\Spec K\to\Spec\Oh_S$ and define
the $p$-adic realizations to be
\[
H^i_p:=H^i_{\et}(\Oh_{S}, j_*T_p).
\]
We will omit the $j_*$, if no confusion is likely. Define 
\[
H_{h,\Z}:=H^m_{\sing}(X\times_{\Q}\C, (2\pi i)^{r-1}\Z)^+
\]
where the $+$ denotes the fixed part under $\Gal(\C/\R)$ 
of the singular cohomology of $X$.
Here $\Gal(\C/\R)$ acts on $\C$ and on $(2\pi i)^{r-1}\Z$.

Finally we need the K-theory of $X$:
Let 
\[
H_{\Mh}=(K_{2r-m-1}(X)\otimes\Q)^{(r)}
\]
be the $r$-th Adams eigenspace of the $2r-m-1$-th Quillen K-theory of 
$X$. 
There are regulator maps due to Beilinson and Soul\'e
\[
r_{\Dh}:H_{\Mh}\otimes_{\Q}\R\to H_{h,\Z}\otimes_{\Z}\R
\]
and
\[
r_p:H_{\Mh}\otimes_{\Q}\Q_p\to H^1_p\otimes_{\Z_p}\Q_p
\]
called the Deligne regulator (see \cite{Be}) and the $p$-adic regulator
(see \cite{So2}).

\begin{rem} Note that because of our assumption $r>\inf(m,\dim(X))$,
the Deligne cohomology coincides with $H_{h,\Z}\otimes_{\Z}\R$
(cf. \cite{Sch} sequence (*) on page 9).
The same condition (together with $m-2r\le -3$) also guarantees that 
$H^i_{\lim}=H^i_p\otimes_{\Z_p}\Q_p$ (cf. \cite{Ka2} 2.2.6 (4)). 
\end{rem}

Let us define local Euler factors for $X$. 
Let for a prime $\pf\nmid p$ in $\Oh_K$
\[
P_{\pf}(V_p, s):={\det}_{\Q_p}(1-\Fr_{\pf}\Norm\pf^{-s}|V_p^{I_{\pf}})
\]
be the characteristic polynomial of the geometric Frobenius $\Fr_{\pf}$ at $\pf$
on the invariants of $V_p$ under the inertia group $I_{\pf}$ at $\pf$.
For $\pf|p$ set
\[
P_{\pf}(V_p, s):={\det}_{\Q_p}(1-\phi_{\pf}^{-1}\Norm\pf^{-s}|D_{\cris}(V_p))
\]
where $D_{\cris}(V_p)):=(V_p\otimes_{\Q_p}B_{\cris})^{\Gal(\Qbar_p/\Q_p)}$
and $\phi_{\pf}$ is the arithmetic Frobenius. Define the 
L-function of $X$ as
\[
L_S(V_p,s):=\prod_{\pf\notin S}P_{\pf}(V_p, s)^{-1}.
\]
Let $V^*_p$ be the dual Galois module of $V_p$.
We now give Kato's formulation of the Tamagawa number
conjecture.
Here and in the rest of the paper the determinants are taken in
the sense of Knudsen and Mumford \cite{Kn-Mu}.
\begin{conj}(cf. \cite{Ka2})\label{iwasawaconj} Let $p\neq 2$ be a prime number, $r,m$ and $S$
be as above. Assume that 
\[
P_{\pf}(V^*_p(1), 0)\neq 0
\]
for all $\pf\in S$ and that $L_S(V^*_p(1),s)$ has an analytic continuation
to all of $\C$. Then:\\
a) The maps $r_{\Dh}$ and $r_p$ are isomorphisms and 
$H^2_p$ is finite.\\
b) $\dim_{\Q}(H_{h,\Z}\otimes\Q)=\ord_{s=0}L_S(V^*_p(1),s)$.\\
c) Let $\eta\in {\det}_{\Z}(H_{h,\Z})$ be a $\Z$-basis and let
$e:=\dim_{\Q}(H_{h,\Z}\otimes\Q)$.
 There is an element $\xi\in {\det}_{\Q}(H_{\Mh})$ such that
\[
r_{\Dh}(\xi)=(\lim_{s\to 0}s^{-e}L_S(V_p^*(1),s))\eta.
\]
This is the ``Beilinson conjecture''.\\
d) Consider $r_p(\xi)\in {\det}_{\Q_p}(H^1_p\otimes_{\Z_p}\Q_p)$,
then $r_p(\xi)$ is a basis of the $\Z_p$-lattice 
\[
{\det}_{\Z_p}(R\Gamma(\Oh_S, T_p))^{-1}\subset
{\det}_{\Q_p}(R\Gamma(\Oh_S, V_p)[-1])\isom {\det}_{\Q_p}(H^1_p\otimes_{\Z_p}\Q_p),
\]
i.e. 
\[
[{\det}_{\Z_p}(H^1_p):r_p(\xi)\Z_p]=\#(H^2_p).
\]
\end{conj}
\begin{rem} a) The assumption in the conjecture is true for abelian varieties
with CM.\\
b) The space $H^0_p$ is zero for weight reasons.\\
c) Part b) follows from the expected shape of the functional equation,
(see e.g. \cite{Sch} proposition page 9).
\end{rem}

As our knowledge of K-theory is limited, let us also
formulate a weak version of  the above conjecture. 
\begin{conj}(weak form of conjecture \ref{iwasawaconj}) There is 
a subspace $H^{\constr}_{\Mh}$ in $ H_{\Mh}$ (the  constructible elements of
$H_{\Mh}$) such that \\
a') $r_{\Dh}$ and $r_p$ restricted  to $H^{\constr}_{\Mh}$ are isomorphisms
and $H^2_p$ is finite.\\
b') same as b)\\
c')  There is an element $\xi\in {\det}_{\Q}(H^{\constr}_{\Mh})$ such that
\[
r_{\Dh}(\xi)=(\lim_{s\to 0}s^{-e}L_S(V^*_p(1),s))\eta.
\]
d') The element $r_p(\xi)$ is a basis of the $\Z_p$-lattice 
\[
{\det}_{\Z_p}(R\Gamma(\Oh_S, T_p))^{-1}\subset
{\det}_{\Q_p}(R\Gamma(\Oh_S, V_p)[-1])\isom {\det}_{\Q_p}(H^1_p\otimes_{\Z_p}\Q_p).
\]
\end{conj}

\subsubsection{Elliptic curves with CM}
Before we formulate the main theorem, we introduce the 
elliptic curves we want to consider.
We follow the notations and conventions in Deninger \cite{Den}.
Let $K$ be an  imaginary quadratic field with ring of integers $\Oh_K$.
Let $E/K$ be an elliptic curve with CM by $\Oh_K$. Note that this implies
that the class number of $K$ is one.
 We fix an isomorphism
\[
\vartheta:\Oh_K\isom \End_{K}(E_K),
\]
such that for $\omega\in \Gamma(E_K, \Omega_{E_K/K})$ and $\alpha\in \Oh_K$
we have $\vartheta^*(\alpha)\omega=\alpha\omega$. 
We fix also an embedding of $K$ into $\C$, such that
the algebraic $j$-invariant of $E$ is the same as the corresponding
complex analytic $j$-invariant of $\Oh_K$. 
Let us denote by 
\[
\psi:\A^*_K\to K^*\subset \C^*
\]
the CM-character or Serre-Tate character of $E_K$ and let $\ff$ be its 
conductor. The elliptic curve $E$ has bad reduction precisely at 
the primes dividing $\ff$. 
Denote by $\bar{\psi}$ the complex conjugate character. Its conductor 
is also $\ff$.
\begin{defn} Fix a prime number $p$. We let  $S$  be the
set of primes in $K$ dividing $p\ff$.
\end{defn}
Associated to $\psi$ is an L-series 
\[
L_S(\psi,s)=\prod_{\pf\nmid p\ff}\frac{1}{1-\frac{\psi(\pf)}{\Norm
\pf^s}}.
\]
We want to relate this to the  L-function $L_S(E,s)$ of 
$E$. Recall the fundamental result of Deuring:
\begin{thm}(see \cite{Si}II 10.5.)
Let $L_S(E/K,s):=L_S(V_p,s)$ be the L-series of the Galois representation
$V_p:=H^1_{\et}(E\times_{K}\Kbar,\Q_p) $ as defined in section
\ref{tamconj}. Then
\[
L_S(E/K,s)=L_S(\psi, s)L_S(\bar{\psi}, s).
\]
\end{thm}
Let $T_pE=\prolim_nE[p^n]$ be the Tate-module of $E$. This is a 
$\Gal(\Kbar/K)$-module.
Then $H^1_{\et}(E\times_{K}\Kbar,\Z_p)\isom\Hom(T_pE,\Z_p)\isom T_pE(-1)$,
where $\Oh_p:=\Oh_K\otimes\Z_p$ acts now conjugate linear on $T_pE$.
There is a canonical isomorphism 
\[
H^1(E\times_{\Q}\C, (2\pi i)^r\Z)^+\isom H^1(E\times_{K}\C, (2\pi i)^r\Z),
\]
where we used the fixed embedding $K\subset \C$.
Let 
\[
H_{\Mh}^i(E,j):=(K(E)_{2j-i}\otimes\Q)^{(j)}
\]
be the $2j-i$-th Quillen K-theory of $E$.

For $m=1$, $r=k+2$ with $k\ge 0$ in the notation of \ref{iwasawaconj}
we have 
\begin{align*}
&H^i_p=H^i_{\et}(\Oh_{S},T_pE(k+1))\\
&H_{h,\Z}=H^1(E\times_{K}\C, (2\pi i)^{k+1}\Z)\\
&H_{\Mh}=H_{\Mh}^2(E,k+2).
\end{align*}
Note that on all these spaces we have a canonical $\Oh_K$-action
and that 
$H_{h,\Z}$ is an  $\Oh_K$-module of rank $1$. It is a result of 
Jannsen (\cite{Ja2} corollary 1) that if $H^2_p$ is finite then 
the free part of $H^1_p$ is an $\Oh_K$-module of rank $1$.

\subsubsection{The main theorem}
Now we can formulate our main result. We let $\Oh_p:=\Oh_K\otimes\Z_p$.
\begin{thm}\label{mainthm} Let $p\neq 2,3$ and
$p\nmid\Norm_{K/Q}\ff$  and $k\ge 0$. Then, there is 
an $\Oh_K$ submodule $\Rh_{\psi}\subset H_{\Mh}$ of rank $1$
such that\\
a) ${\det}_{\Oh_K}(r_{\Dh}(\Rh_{\psi}))\isom L_S^*(\bar{\psi},-k){\det}_{\Oh_K}(H_{h,\Z})$ in  ${\det}_{\Oh_K\otimes\R}(H_{h,\Z}\otimes\R)$\\
and \\
b) The map $r_p$ induces an isomorphism
\[
{\det}_{\Oh_p}(\Rh_{\psi}) \isom{\det}_{\Oh_p}(R\Gamma_{\et}(\Oh_{S},T_pE(k+1)))^{-1}.
\]
Here $L^*(\bar{\psi},-k)=\lim_{s\to -k}\frac{L(\bar{\psi},s)}{s}$ denotes the 
leading coefficient of the Taylor series of $L(\bar{\psi},s)$ at $-k$.
Moreover, if $H^2_p$ is finite, $r_p$ is injective on $\Rh_{\psi}$ and
\[
{\det}_{\Oh_p}H^1_p/r_p(\Rh_{\psi})\isom {\det}_{\Oh_p}H^2_p.
\]
\end{thm}
\begin{rem} i) Part a) was proven by Deninger in \cite{Den}. This is
the Beilinson conjecture for Hecke characters.\\
ii) For $k=0$ and CM elliptic curves defined over $\Q$ with CM by $\Oh_K$
and $p$ regular, part b) was proven in \cite{Bl-Ka}. They
used an ad hoc computation of the $p$-adic realization of 
the K-theory elements, which does not generalize.\\
iii) Let us explain part b) in more detail: The Soul\'e regulator 
\[
r_p:\Rh_{\psi}\to H^1_p\otimes\Q_p
\]
extends to a map to $R\Gamma_{\et}(\Oh_{S},T_pE(k+1))[-1]\otimes\Q_p$
because $H^0_p$ is zero for weight reasons. The determinant of this
complex has as $\Oh_p$-lattice the determinant of 
$R\Gamma_{\et}(\Oh_{S},T_pE(k+1))[-1]$, which is 
${\det}_{\Oh_p}(R\Gamma_{\et}(\Oh_{S},T_pE(k+1)))^{-1}$.
Part b) means that the determinant of the complex 
\[
\Rh_{\psi}\xrightarrow{r_p}R\Gamma_{\et}(\Oh_{S},T_pE(k+1))[-1]
\]
is trivial. If $H^2_p$ is finite, and $r_p\otimes\Q_p$ an isomorphism,
we get from this ${\det}_{\Oh_p}(H^1_p/r_{p}(\Rh_{\psi}))\isom {\det}_{\Oh_p}H^2_p$.\\
iv) Note that  $H^2_p$ is finite for almost all $k\ge 0$
or for $p$ regular.
See the next section for remarks  on the finiteness of $H^2_p$.\\
v) From this result for the set of primes $S$ we get 
it for all other set of primes which contain $S$.
This follows from 
\cite{Ka2} 4.11. 
\end{rem}

As a corollary
we get a result about  $\Z_p$ determinants and 
the L-function of $E/K$:
\begin{cor}
Under the conditions of the theorem \\
a) ${\det}_{\Z}(r_{\Dh}(\Rh_{\psi})) =L_S^*(E/K,-k){\det}_{\Z}(H_{h,\Z})$\\
and that\\
b) ${\det}_{\Z_p}(r_{p}(\Rh_{\psi})) ={\det}_{\Z_p}(R\Gamma_{\et}(\Oh_{S},T_pE(k+1)))^{-1}.$\\
Here $L^*(E/K,-k)=\lim_{s\to -k}\frac{L(E/K,s)}{s^2}$ denotes the 
leading coefficient of the Taylor series of $L(E/K,s)$ at $-k$.
\end{cor}
\bew This follows from the theorem and the remark that if we multiply 
an $\Oh_K$-module with an element $L^*(\psibar, -k)$
from $\Oh_K\otimes\R$, then the
determinant is multiplied by the norm $\Norm_{\Oh_K\otimes\R/\R}(L^*(\psibar, -k))$.
But $\Norm_{\Oh_K\otimes\R/\R}(L^*(\psibar, -k))=L^*(E/K,-k)$.
Part b) is obvious.
\bewende

Here is a short overview of the proof: We start by recalling 
Deninger's definition of K-theory elements and his main
result about the relation of these to the $L$-value.
Then we use an idea of Soul\'e to construct a  submodule of 
\[
H^1_{\et}(\Oh_{S},T_pE(k+1))
\]
 via elliptic units.
Iwasawa theory allows us to compute the index of this
submodule. The proof concludes with the comparison of this submodule
and $\Rh_{\psi}$.
This is the main step in the proof which needs the theory
developed in part \ref{seccomp} and in particular the explicit description
of the elliptic polylogarithm of theorem  \ref{polastheta}.
The injectivity of $r_p$ in the case that $H^2_p$ is finite, will
be proved in section \ref{h2finite}.
\subsubsection{Some remarks concerning the finiteness of $H^2$}
Note that we do {\em not} prove that $H^2_p$ is finite. 
Nevertheless, there is the following result of Soul\'e:
\begin{thm}(\cite{So} 1.5 proposition 3) For fixed $p$ the 
group 
\[
H^2(\Oh_{S},T_pE(k+1))
\]
 is finite for almost all $k$. 
\end{thm}
For regular $p$, we have results of Soul\'e  and Wingberg:
\begin{thm}(\cite{So2} 3.3.2, \cite{Win} cor. 2) Let $p$ be a regular
prime for $E$ (see e.g.\cite{So2} 3.3.1 for the definition of regular), then
\[
H^2_{\et}(\Oh_S, E[p^{\infty}](k+1))=0.
\]
\end{thm}
It is easy to see that this vanishing implies the finiteness of $H^2_p$ 
(cf. \cite{Ja2} lemma 1). \\
\begin{rem} In \cite{Ja2} it is conjectured that $H^2_p$ is always finite.
\end{rem}

\subsection{Review of the Deninger elements for CM elliptic curves over an
imaginary quadratic field}
Here we  describe briefly the construction by Deninger \cite{Den} of 
the elements in K-theory, which interpret
the L-value up to rational numbers as predicted by Beilinson's conjecture.

We are interested in the L-values $L(\bar{\psi}, k+2)$ with $k\ge 0$.
Thus according to \ref{iwasawaconj} we need an element in 
$H^{2}_{\Mh}(E,k+2)$. 

\subsubsection{The Deninger--Beilinson construction}
We fix an algebraic differential $\omega\in H^0(E,\Omega_{E/K})$ and
let $\Gamma$ be its period lattice. Then we have an isomorphism
\begin{align*}
E(\C)&\to \C/\Gamma\\
z&\mapsto \int_0^z\omega
\end{align*}
using the fixed embedding $K\subset \C$. This isomorphism is equivariant
for the action of complex multiplication and because $j(E)=j(\Oh_K)$
the lattice is of the form $\Gamma=\Omega\Oh_K$ for some $\Omega\in\C^*$.
Fix an $\Oh_K$ generator $\gamma\in H_1(E(\C),\Z)$, then
\[
\Omega=\int_{\gamma}\omega.
\]
Recall that $\ff$ is the conductor of $\psi$ and the locus of bad
reduction of $E$. 
Let $\Z[E[\ff]\ohne 0]$ be the  group of divisors with support
in the $\ff$-torsion points without $0$ of $E$, defined 
over $K$. Beilinson defines a map:
\begin{thm}(\cite{Be1} There is a non-zero map, a variant of 
the Eisenstein symbol,
\[
\Z[E[\ff]\ohne 0]\xrightarrow{\Eh^{2k+1}_{\Mh}}H^{2k+2}_{\Mh}(E^{2k+1},2k+2),
\]
where $E^n:=E\times_K\ldots \times_KE$.
\end{thm}
Deninger constructs a projector
\[
\Kh_{\Mh}:H^{2k+2}_{\Mh}(E^{2k+1},2k+2)\to H^{2}_{\Mh}(E,k+2)
\]
 as follows: Let $d_K$ be the discriminant of $K$ and $\sqrt{d_K}$
be a square root of $d_K$. Complex multiplication gives a map
\[
\delta=(\id,\vartheta(\sqrt{d_K})):E\to E\times_KE
\]
and taking this $k$-times gives
$\delta^k\times\id:E^k\times_KE\to E^{2k}\times_KE$. Then
\[
\Kh_{\Mh}=\pr_*\verk(\delta^k\times\id)^*
\]
where $\pr$ is the projection $E^k\times_KE\to E$ onto the last component.
Hence we get a map 
\[
\Kh_{\Mh}\verk\Eh^{2k+1}_{\Mh}:\Z[E[\ff]\ohne 0]\to  H^{2}_{\Mh}(E,k+2).
\]

\subsubsection{The Beilinson conjecture for CM elliptic curves}\label{Beilinson}
Following Deninger we define an element $\beta$ in $\Z[E[\ff]\ohne 0]$. 
Let $K(\ff)$ be the ray class field associated to $\ff$ and
note that $K(\ff)=K(E[\ff])$.

Let $f$ be a generator of $\ff$. Then
\begin{equation}\label{tdefn}
\Omega f^{-1}\in  \ff^{-1}\Gamma
\end{equation}
 defines an element in $E[\ff](K(\ff))$. 
This gives a divisor $(\Omega f^{-1})$ in 
$\Z[E[\ff]\ohne 0]$ defined over $K(\ff)$ on which the Galois
group $\Gal(K(\ff)/K)$ acts. We  define:
\[
\beta:=\Norm_{K(\ff)/K}((\Omega f^{-1})).
\]
This is a divisor defined over $K$.
Recall that $\gamma$ is an $\Oh_K$ generator of $H_1(E(\C),\Z)$.
By Poincar\'e duality we have an isomorphism (conjugate linear for the $\Oh_K$-action)
\[
H^1(E(\C),\Z(k+1))\isom\Hom(H^1(E(\C),\Z), \Z(k))=H_1(E(\C),\Z(k)).
\]
Denote by $\eta$ the $\Oh_K$ generator of $H^1(E(\C),\Z(k+1))$ corresponding
to $(2\pi i)^k\gamma$ under this isomorphism.
We can now formulate the main result of \cite{Den} in our case:
\begin{thm}(\cite{Den} thm. 11.3.2)\label{deligneregulator}
Let $\beta$ and $\eta$ be as above and define
\[
\xi:=(-1)^{k-1}\frac{(2k+1)!}{2^{k-1}}\frac{L_p(\bar{\psi},-k)^{-1} }{\psi(f)\Norm_{K/\Q}\ff^k}
\Kh_{\Mh}\verk\Eh^{2k+1}_{\Mh}(\beta)\in
 H^{2}_{\Mh}(E,k+2),
\]
where $L_p(\bar{\psi},-k)$ is the Euler factor of $\bar{\psi}$ at $p$, 
evaluated at $-k$.
Then
\[
r_{\Dh}(\xi)=L^*_S(\psibar, -k)\eta\in H^1(E\times_K\C,(2\pi i)^{k+1}\R),
\]
where $L^*_S(\psibar, -k)=\lim_{s\to -k}\frac{L_S(\psibar, s)}{s}$.
\end{thm}
Note that $L_{\pf}(\bar{\psi},-k)=1$, if $\pf|\ff$.
We can now define the space $\Rh_{\psi}$ of the main theorem
\ref{mainthm}
\begin{defn}\label{Rhdefn}
We define
\[
\Rh_{\psi}:=\xi\Oh_K\subset  H^{2}_{\Mh}(E,k+2)
\]
to be the $\Oh_K$-submodule of $  H^{2}_{\Mh}(E,k+2)$ generated by 
$\xi$.
\end{defn}
Note that by the above theorem, $\Rh_{\psi}$ is an $\Oh_K$-modules of 
rank $1$.
\begin{cor} With the above notation
\[
r_{\Dh}({\det}_{\Z}(\Rh_{\psi}))=L^*_S(E/K,-k){\det}_{\Z}(H^1(E(\C),\Z(k+1))).
\]
where $S$ is the set of primes in $K$ dividing $p\ff$.
\end{cor}
\bew This follows from the theorem and the remark that if we multiply 
an $\Oh_K$-module with an element $L^*_S(\psibar, -k)$
from $\Oh_K\otimes\R$, then the
determinant is multiplied by the norm $\Norm_{\Oh_K\otimes\R/\R}(L^*_S(\psibar, -k))$.
But $\Norm_{\Oh_K\otimes\R/\R}(L^*_S(\psibar, -k))=L^*_S(E/K,-k)$.
\bewende

\subsubsection{The space $r_p(\Rh_{\psi})$ in terms of the specialization of 
the elliptic polylog}
Recall from definition \ref{Rhdefn} that the space $\Rh_{\psi}$
is generated as an $\Oh_K$-module by the element $\xi$ 
from theorem
\ref{deligneregulator}. The element $\xi$ is up to some factors
of the form $\Kh_{\Mh}\verk\Eh^{2k+1}_{\Mh}(\beta)$.
Let us define 
\[
t:=\Omega f^{-1}
\]
with the notation from \ref{tdefn}.
This is an $\Norm_{K/\Q}\ff$-division point. Then
we have $\beta=\Norm_{K(\ff)/K}((t))$.
Now let 
\[
(\beta^*\pol_{\Q_p})^{2k+1}\in H^1_{\et}(\Oh_S, \Sym^{2k+1}\shH_{\Q_p}(1))
\]
be the specialization of the polylogarithm as defined 
in \ref{polylogspecialization}.
We have the following comparison theorem:
\begin{thm}\label{eiseq} There is an equality
\[
r_p(\Eh^{2k+1}_{\Mh}(\beta))=-\Norm\ff^{4k+2}(\beta^*\pol_{\Q_p})^{2k+1}
\]
in $H^1_{\et}(\Oh_S, \Sym^{2k+1}\shH_{\Q_p}(1))$.
\end{thm}
\bew The formula is the combination of 
two results: Theorem 2.2.4 in \cite{Hu-Ki}, which states that
\[
r_p(\Eh is^{2k+1}_{\Mh}(\rho\beta))=-\Norm\ff^{2k}(\beta^*\pol_{\Q_p})^{2k+1}
\]
where $\Eh is^{2k+1}_{\Mh}$ is Beilinson's Eisenstein symbol
and $\rho$ the horospherical map. 
Note that what is here called $(\beta^*\pol_{\Q_p})^{2k+1}$ is in loc. cit.
$(\beta^*\pol_{\Q_p})^{2k+2}$.
Furthermore, according to \cite{Den1} formula
3.35., 
\[
\Eh^{2k+1}_{\Mh}(\beta)=\Norm\ff^{2k+2}\Eh is^{2k+1}_{\Mh}(\rho\beta).
\]
Note that the formula 3.35. in \cite{Den1} uses an other normalization
of the horospherical map and that there is a factor $\Norm\ff$ missing because 
of a wrong normalization of the residue map (the residue of $\frac{dq}{q}$
in formula $(3.7)$ in loc. cit. is not $1$ but $N$).
\bewende

We have now two tasks: To compute the specialization of the
elliptic polylog and to identify the \'etale cohomology groups
to compute the ``index'' of $r_p(\Rh_{\psi})$ in 
$H^1_{\et}(\Oh_S, \Sym^{2k+1}\shH_{\Z_p}(1))$. The answer to these
problems involves elliptic units and we will start to use Iwasawa
theory and Rubin's proof of the main conjecture to identify the 
the \'etale cohomology groups. 

\section{Iwasawa theory}
This section treats the relation between certain Iwasawa modules
and \'etale cohomology. Rubin's ``main conjecture'' is used in an
essential way. We first review the results of Rubin and
then use an idea of Soul\'e to produce elements in \'etale
cohomology using elliptic units. An idea of Kato, which 
partly goes back to Soul\'e as well, allows
to compare the elliptic units and \'etale cohomology. The
main result of this part is theorem \ref{cohdescr}.

\subsection{Review of the ``main conjecture'' of Iwasawa theory
for CM elliptic curves}\label{secmainconj}
In this section we review the  ``main conjecture'' of Iwasawa theory
for CM elliptic curves proved by Rubin \cite{Ru1}. This will be
used in section \ref{secreduction}  to reduce the
Tamagawa number conjecture to an ``index computation''.  

\subsubsection{Definition of the Iwasawa modules}
We  follow Rubin \cite{Ru1}:
Let $E/K$ be as before an elliptic curve with CM by $\Oh_K$,
$K$ an imaginary quadratic field. We fix an embedding
of $K$ into $\C$ and view $K$ as a subfield of $\C$. 
Fix a prime $p\nmid \#\Oh_K^*$ and
a prime $\pf$ of $\Oh_k$ lying over $p$ and  denote by $E[p^n]$ the 
$p^n$-torsion points of $E$. Let $K_n:=K(E[p^{n+1}])$ be the 
extension field defined by these torsion points and $\Kinfty:=\indlim_nK_n$.
Denote the ring of integers in these fields by $\Oh_n$ (resp. $\Ohinfty$).
Then $\Delta:=\Gal(K_0/K)$
has order prime to $p$ and $\Gamma:=\Gal(\Kinfty/K_0)$ is isomorphic to
$\Z_p^2$. Let $\Gh:=\Gal(\Kinfty/K)$ be the 
Galois group of the extension $\Kinfty/K$. Then $\Gh\isom\Delta\times\Gamma$.
Define $\Ah_n$ to be the $p$-part of the 
ideal class group of $K_n$, $\Eh_n$ to be the group of 
global units $\Oh_n^*$ of $K_n$ and $\Uh_n^{\pf}$ 
 the local units of $K_n\otimes_KK_{\pf}$
which are congruent to $1$ modulo the primes above $\pf$. 
For every
prime $v$ of $K_n$ above $\pf$ there
is an exact sequence
\begin{equation}\label{localunits}
1\to \Uh_{n,v}\to K_{n,v}^*\to \Z\times \kappa_n^*\to 1
\end{equation}
and $\Uh_n^{\pf}=\bigoplus_{v|\pf}\Uh_{n,v}$.
Here $\Uh_{n,v}$ are the local units congruent to 
$1$ modulo $v$ and $\kappa_n$ is the residue class field of $K_{n,v}$.
Let $\Ch_n$ be the elliptic units in $K_n$ as defined in \cite{Ru1} paragraph 1.
We recall their definition. For every ideal $\af\in\Oh_K$ prime to $6$
consider  the function $\theta_{\af}(z)$ that will be
defined in  \ref{thetadefn}.
 The function  $\theta_{\af}(z)$
is a $12$-th root of the function in \cite{deSh} II.2.4.

Let $t:=\Omega f^{-1}$  and $\af$ be 
an ideal prime to $6\ff$. 
\begin{defn}\label{ellipticunits}(cf.\cite{Ru3} 11.2)
Let  $C_n$ be the subgroup of units generated over 
$\Z[\Gal(K_n/K)]$ by
\[
\prod_{\sigma\in \Gal(K(\ff)/K)}\theta_{\af}(t^{\sigma}+h_n),
\]
where $\af$ runs through
all ideals prime to $6p\ff$, 
$K(\ff)$ is the ray class field defined by $\ff$ and 
$h_n$ is a primitive $p^n$-torsion point.
Define
\[
\Ch_n:=\mu_{\infty}(K_n)C_n,
\]
the group of {\em elliptic units} of $K_n$.
\end{defn}
Denote by $\Ehbar_n$ and $\Chbar_n$ the closures of 
$\Eh_n\cap\Uh_n^{\pf}$ resp. $\Ch_n\cap\Uh_n^{\pf}$ in $\Uh_n^{\pf}$. Finally define
\[
\Ahinfty:=\prolim_n\Ah_n,\;\;\; \Ehbarinfty:=\prolim_n\Ehbar_n,\;\;\;
 \Chbarinfty:=\prolim_n\Chbar_n,\;\;\;
 \Uhinfty^{\pf}:=\prolim_n\Uh_n^{\pf}
\]
where the limits are taken with respect to the norm maps. 
Denote by $\Minfty^{\pf}$ the maximal abelian $p$-extension of $\Kinfty$
which is unramified outside of the primes above $\pf$, and write
$\Xhinfty^{\pf}:=\Gal(\Minfty^{\pf}/\Kinfty)$. Global class field theory 
gives an exact sequence
\begin{equation}\label{iwasawaseq}
0\to \Ehbarinfty/ \Chbarinfty\to  \Uhinfty^{\pf}/ \Chbarinfty\to \Xhinfty^{\pf}\to
\Ahinfty\to 0.
\end{equation}
Define the Iwasawa algebra 
\[
\Z_p[[\Gh]]:=\prolim_n\Z_p[[\Gal(K_n/K)]]
\]
this has an action of $\Z_p[\Delta]$. For any irreducible $\Z_p$-representation
$\chi$ of $\Delta$, let 
\[
e_{\chi}:=\frac{1}{\#\Delta}\sum_{\tau\in\Delta}\Tr(\chi(\tau))\tau^{-1}
\in\Z_p[\Delta]
\]
and let for every $\Z_p[\Delta]$-module $Y$ be $Y^{\chi}:=e_{\chi}Y$ 
be the $\chi$-isotypical component.
In particular we define 
\[
\Lambda^{\chi}:=\Z_p[[\Gh]]^{\chi}=R_{\chi}[[\Gamma]]
\]
where $R_{\chi}$ is the ring of integers in the unramified extension
of $\Z_p$ of  degree $\dim(\chi)$. As we will work with
$\Z_p[[\Gamma]]\otimes \Oh_p$-modules,  we let
\[
\Lambda:=\Oh_p[[\Gamma]].
\]
Then $ \Ehbarinfty^{\chi}$, $\Chbarinfty^{\chi}$, ${\Uhinfty^{\pf}}^{\chi}$,
$\Ahinfty^{\chi}$ and $\Xhinfty^{\chi}$ are 
finitely generated $\Lambda^{\chi}$-modules (see \cite{Ru1} paragraph 5).
The modules $\Ahinfty^{\chi}$ and $ \Ehbarinfty^{\chi}/\Chbarinfty^{\chi}$
are even torsion $\Lambda^{\chi}$-modules.

\subsubsection{Rubin's ``main conjecture'' for imaginary quadratic fields}
We have the following lemma due to Kato:
\begin{lemma}(see \cite{Ka2} proposition 6.1.)
Let $Y$ be a finitely generated torsion $\Lambda^{\chi}$-module.
Then 
\[
{\det}_{\Lambda^{\chi}}(Y)=\car (Y),
\]
where $\car(Y)$ is the usual characteristic ideal in Iwasawa theory
(see e.g. \cite{Ru1} paragraph 4) and the determinant is taken
in the sense of \cite{Kn-Mu}.
\end{lemma}
With this lemma we can formulate the main result of \cite{Ru1} 
as follows:
\begin{thm}(\cite{Ru1} theorem 4.1.)\label{mainconj} Let $p\nmid \#\Oh_K^*$.\\
i) Suppose that $p$ splits in $K$, then
\[
{\det}_{\Lambda^{\chi}}(\Ahinfty^{\chi})=
{\det}_{\Lambda^{\chi}}(\Ehbarinfty^{\chi}/\Chbarinfty^{\chi}).
\]
ii) Suppose that $p$ remains prime or ramifies in $K$ and that
$\chi$ is nontrivial on the decomposition group of $\pf$ in $\Delta$, then
\[
{\det}_{\Lambda^{\chi}}(\Ahinfty^{\chi})=
{\det}_{\Lambda^{\chi}}(\Ehbarinfty^{\chi}/\Chbarinfty^{\chi}).
\]
\end{thm}
Using the theory of the determinant and the exact sequence
(\ref{iwasawaseq}) we get:
\begin{cor}
In the situation of the theorem \ref{mainconj},
\[
{\det}_{\Lambda^{\chi}}({\Xhinfty^{\pf}}^{\chi})=
{\det}_{\Lambda^{\chi}}({\Uhinfty^{\pf}}^{\chi}/\Chbarinfty^{\chi}).
\]
\end{cor}
We need a variant of this. 
Let ${\Xhinfty}$ be the Galois group of the maximal abelian
$p$-extension $\Minfty^p$ of $\Kinfty$ which is unramified outside of
the primes above $p$. In the case where $p$ is inert or ramified
this is the same as $\Xhinfty^{\pf}$. Define also
\[
\Uhinfty:=\Uhinfty^{\pf}\times\Uhinfty^{\pf^*}
\]
if $p=\pf\pf^*$ is split, and 
\[
\Uhinfty:= \Uhinfty^{\pf}
\]
if $p$ is inert or ramified. Let similarly
$\Yh_n$ be the $p$-adic completion of 
$(K_n\otimes\Q_p)^*$ and $\Yhinfty:=\prolim_n\Yh_n$.
We have an inclusion 
$\Uhinfty\subset \Yhinfty$.
Class field theory gives
\begin{equation}\label{iwasawaseq2}
0\to \Ehbarinfty/ \Chbarinfty\to \Uhinfty
/ \Chbarinfty\to \Xhinfty\to
\Ahinfty\to 0.
\end{equation}
where $\Chbarinfty$ is diagonally embedded into 
$\Uhinfty^{\pf}\times\Uhinfty^{\pf^*}$ if $p$ is split. 
On this sequence acts $\Gh$ and we get:
\begin{cor}\label{mainconjcor}
\[
{\det}_{\Lambda^{\chi}}({\Xhinfty}^{\chi})=
{\det}_{\Lambda^{\chi}}({\Uhinfty}^{\chi}/\Chbarinfty^{\chi}).
\]
\end{cor}

\begin{lemma}\label{UYrel} Let $p\nmid\Norm\ff$ be a prime. If
$p$ splits in $K$, the inclusion $\Uhinfty\to \Yhinfty$
is an isomorphism and if $p$ is inert or ramified in $K$,
there is an exact sequence
\[
0\to \Uhinfty\to \Yhinfty \to \Z_p[\Delta/\Delta_p]\to 0,
\]
where $\Delta_p$ is the decomposition group of $p$ in $\Delta=\Gal(K_0/K)$.
\end{lemma}
\bew
We have exact sequences
\[
1\to \Uh_{n,v}\to K_{n.v}^*\to \Z\times \kappa_n^*\to 1
\]
where $\kappa_n$ is the residue class field of $ K_{n.v}$.
By definition $\Uh_n=\bigoplus_{v|p}\Uh_{n,v}$.

As the order of the residue class field $\kappa_n^*$ is prime
to $p$, we have an exact sequence 
\[
0\to \prolim_n\Uh_{n,v}/p^n\to \prolim_n K_{n,v}^*/p^n\to \Z_p\to 0.
\]
As $E$ has good reduction at $p$, we now how $p$ decomposes
in $K_n$ (see \cite{Ru2} prop. 3.6). 
If $p$ is split, the ramification degree of $v$ in $K_{n+1}$ is $p$
 and the degree of $K_{n+1}$ over $K_n$ is $p^2$.
Hence the norm map induces multiplication by $p$ on $\Z_p$
and the inverse limit over these maps is zero.
This gives the first claim. 
In the case where $p$ is inert of ramified, $v$ is totally ramified in 
$K_{n+1}$ and 
the norm map from $K_{n',v}^*\to K_{n,v}^*$ induces the
identity on $\Z_p$.
Putting these sequences together for all $v|p$ and using 
$\bigoplus_{v|p}\Z_p=\Z_p[\Delta/\Delta_p]$ gives the result.
\bewende

\subsection{Reductions via Iwasawa theory}\label{secreduction}

In this section we use Rubin's  ``main conjecture'' of Iwasawa theory 
to reduce the Bloch-Kato conjecture to a comparison between 
the space $\Rh_{\psi}$ of \ref{Rhdefn} and the elliptic units $\Chbarinfty$.

We have a subspace
\[
r_p(\Rh_{\psi})\subset  H^1_{\et}(\Oh_S,V_p) 
\]
where $\Oh_S=\Oh_K[\frac{1}{S}]$ 
and $V_p=T_pE(k+1)\otimes\Q_p$. Recall that $S$ is the set of
primes of $K$ dividing $p\ff$. We want to 
compute the relation  of the submodule
$r_p(\Rh_{\psi})$ to $H^1_{\et}(\Oh_S,T_pE(k+1))$.
Our method, which is inspired by Kato's paper \cite{Ka2},
relates the module $r_p(\Rh_{\psi})$ to a certain submodule defined
by the elliptic units $\Chbarinfty$. This submodule in turn is defined using
an idea of Soul\'e.
Our aim is to relate the determinant of $\Chbarinfty\otimes T_pE(k)$ to the
determinant of $R\Gamma(\Oh_S,T_pE(k+1))$ (see theorem \ref{cohdescr} for the
exact formulation).

In this section $p$ is always a prime which does not divide $\#\Oh_K^*$
and where $E$ has good reduction over the primes above $p$, i.e. 
$p\nmid\Norm\ff$.

Denote by abuse of notation by $S_p$ the
set of primes over $p$ in the 
ring $\Oh_n$ for every $n$ and by
$\Oh_{n,S_p}$ the ring of integers in $K_n$ where
the primes above $p$ are inverted. We define $\Oh_{\infty,S_p}:=\indlim_n
\Oh_{n,S_p}$. Similarly we define $\Oh_{n,S}$ and $\Oh_{\infty,S}$.

\subsubsection{Review of the Soul\'e elements}
 We keep the notations from 
the section \ref{secmainconj}.

Denote by $T_{p}E=\prolim_nE[p^n]$
the Tate module of $E$ and let $T_pE(k):=T_pE\otimes\Z_p(k)$ its Tate twist.
This is a $\Oh_p\otimes\Z_p[[\Gh]]$-module.
We start by defining a map in the spirit of Soul\'e
\[
\Chbarinfty\otimes_{\Z_p} T_{p}E(k)\to H^1(\Oh_{S},T_{p}E(k+1))
\]
here $T_{p}E(k+1)$ is a sheaf on $\Oh_S$ because it is unramified outside of
$S$.
Write 
\[
 H^1(\Oh_{S},T_{p}E(k+1))=\prolim_rH^1(\Oh_{S},E[p^{r+1}](k+1))
\]
and  let for a norm compatible system of elliptic units
$(\theta_r)_r$ and an element $(t_r)_r\in T_{p}E(k+1) $
\[
e_p((\theta_r\otimes t_r)_r):=(\Norm_{K_r/K}(\theta_r\otimes t_r))_r
\]
where $\theta_r\otimes t_r$ is an element in
\[
\Oh_{r,S}^*/(\Oh_{r,S}^*)^{p^{r+1}}\otimes
E[p^{r+1}](k)\subset H^1(\Oh_{r,S},E[p^{r+1}](k+1))
\]
(the inclusion comes from Kummer theory) and $\Norm_{K_r/K}$ is
the norm map on the cohomology.
According to Soul\'e (\cite{So2} lemma 1.4) this gives a projective system of
 elements in $ H^1(\Oh_{S},T_{p}E(k+1))$. 
The map $e_p$ factors through 
the coinvariants under $\Gh$, so that we can make the
following definition:
\begin{defn}\label{souleelements}
The {\em Soul\'e elliptic elements} are defined by the map
\[
e_{p}:(\Chbarinfty\otimes T_{p}E(k))_{\Gh}
\to H^1(\Oh_{S},T_{p}E(k+1)).
\]
\end{defn}
We want to investigate to what extend this map gives generators
 for $ H^1(\Oh_{S},T_{p}E(k+1))$, i.e. what is the kernel and cokernel
of the map $e_p$. 

\subsubsection{The Tate-Poitou localization sequence}
Our main tool in describing $ H^1(\Oh_{S},T_{p}E(k+1))$ in terms
of elliptic units will be the Tate-Poitou localization sequence. 
It is convenient for us to write down a derived category version of it.
 
For technical reasons we have to
work first over $K_0$.
The reason for  this is that over $K_0$ the module $T_{p}E(k+1)$ is unramified 
outside of the primes above $p$:
\begin{lemma}(\cite{Ru2}1.3)\label{redlemma} If $p\nmid \#\Oh_K^*$, then over $K_0$ the 
elliptic curve $E$ has good
reduction at all places not dividing $p$. In particular there
exists a model of $E$ over $\Oh_{0,S_p}$ and $T_{p}E(k+1)$ is unramified.
\end{lemma}
The localization sequence now reads as follows (see \cite{Ka2} (6.3)).
 Here $^*$ is the Pontryagin dual 
$\Hom_{\Oh_p}(\_,\Q_p/\Z_p\otimes_{\Z_p}\Oh_p)$.
\begin{align*}
R\Gamma_{\et}(\Oh_{0,S_p}, T_{p}E(k+1))\to&
R\Gamma_{\et}(K_0\otimes
 \Q_{p}, E[p^{\infty}](-k))^*[-2]\\
&\to R\Gamma_{\et}(\Oh_{0,S_p}, E[p^{\infty}](-k))^*[-2]\to,
\end{align*}
where we have used the identification 
\[
E[p^{\infty}](-k)=\Hom_{\Oh_p}( T_{p}E(k+1),\Q_p/\Z_p(1)\otimes\Oh_p).
\]
Our next task is to rewrite this Tate-Poitou sequence in terms of 
Iwasawa theory.

\subsubsection{Identification of some Galois cohomology groups with
Iwasawa modules}
 Let us define
\begin{align*}
H^1_{\et}(\Kinfty\otimes
 \Q_{p}, E[p^{\infty}](-k)):&=\indlim_nH^1_{\et}(K_n\otimes
\Q_{p}, E[p^{\infty}](-k))\\
&=\indlim_n\bigoplus_{v|p}H^1_{\et}(K_{n,v}, 
E[p^{\infty}](-k)).
\end{align*}
Note that there are only finitely many primes above $p$ in $\Kinfty$.

\begin{prop}\label{iwasawaident}
There are isomorphisms of $\Oh_p[[\Gh]]$-modules
\begin{align*}
\Xhinfty\otimes_{\Z_p} T_{p}E(k)&\isom H^1_{\et}(\Oh_{\infty,S_p}, E[p^{\infty}](-k))^*\\
\Yhinfty\otimes_{\Z_p} T_{p}E(k)&\isom H^1_{\et}(\Kinfty\otimes
 \Q_{p}, E[p^{\infty}](-k))^*,
\end{align*}
where $^*$ 
is the Pontryagin dual $\Hom_{\Oh_p}(\_,\Q_p/\Z_p\otimes_{\Z_p}\Oh_p)$.
\end{prop}
\bew 
We have
\begin{align*}
H^1_{\et}(\Oh_{\infty,S_p}, E[p^{\infty}](-k))^*&
=\Hom (\Gal(\Kbar/\Kinfty), E[p^{\infty}](-k))^*\\
&=\Hom (\Gal(\Minfty^p/\Kinfty), E[p^{\infty}](-k))^*\\
&=\Xhinfty\otimes T_{p}E(k).
\end{align*}
In the local case, we have an isomorphism
\begin{align*}
 H^1_{\et}(K_n\otimes
 \Q_{p}, E[p^{n}](-k))&=\bigoplus_{v|p}H^1_{\et}(K_{n,v}, E[p^{n}](-k))\\
&=\bigoplus_{v|p}\Hom (\Gal(\bar{K}_{n,v}/K)^{\ab}, E[p^{n}](-k))
\end{align*}
By class field  theory
\[
\Hom (\Gal(\bar{K}_{n,v}/K)^{\ab}, E[p^{n}](-k))^*\isom K_{n,v}^*/p^n\otimes
E[p^{n}](k)
\]
so that
\[
 H^1_{\et}(\Kinfty\otimes
 \Q_{p}, E[p^{\infty}](-k))^*=\bigoplus_{v|p}\Yhinfty\otimes T_{p}(k).
\]
\bewende

\subsubsection{Rewriting the Tate-Poitou localization sequence in terms
of Iwasawa theory}

To proceed further, we need the following vanishing result.
\begin{prop}\label{h2van}
The groups 
\[
H^2_{\et}(\Kinfty\otimes
 \Q_{p}, E[p^{\infty}](-k))\mbox{\quad and \quad}H^2_{\et}(\Oh_{\infty,S_p}, E[p^{\infty}](-k))
\]
are zero.
\end{prop}
\bew By local duality we have
\[
H^2_{\et}(K_n\otimes \Q_{p}, E[p^{\infty}](-k))^*\isom
H^0_{\et}(K_n\otimes \Q_{p}, T_{p}E(k+1))=0.
\]
On the other hand it is a result of Schneider \cite{Sch1} 4.1 that the 
cohomology group $H^2_{\et}(\Oh_{\infty,S_p},\Q_p/\Z_p(-k))$ is zero.
As 
\[
H^2_{\et}(\Oh_{\infty,S_p}, E[p^{\infty}](-k))=
H^2_{\et}(\Oh_{\infty,S_p},\Q_p/\Z_p(-k))\otimes_{\Z_p} T_pE
\]
this proves our claim.
\bewende
This vanishing result implies that we get actually a map from the 
Iwasawa modules to complexes computing the Galois cohomology.
\begin{cor}\label{extri}
There are exact triangles
\begin{align*}
\Yhinfty\otimes T_{p}E(k)[1]&\to R\Gamma_{\et}(\Kinfty\otimes
 \Q_{p}, E[p^{\infty}](-k))^*\to H^0_{\et}(\Kinfty\otimes
 \Q_{p}, E[p^{\infty}](-k))^*\\
\Xhinfty\otimes T_{p}E(k)[1]&\to R\Gamma_{\et}(\Oh_{\infty,S_p}, E[p^{\infty}](-k))^*\to H^0_{\et}(\Oh_{\infty,S_p}, E[p^{\infty}](-k))^*
\end{align*}
\end{cor}
\bew The propositions \ref{iwasawaident} and \ref{h2van} show that
we have a canonical map 
from $\Yhinfty\otimes T_{p}E(k)[1]$ to  $ R\Gamma_{\et}(\Kinfty\otimes
  \Q_{p}, E[p^{\infty}](-k))^*$ 
because the second cohomology vanishes.
The same argument gives the result for $\Xhinfty\otimes T_{p}E(k)[1]$.
\bewende

To relate these groups to the cohomology groups of $\Oh_{S_p}$
we want to take the coinvariants under $\Gamma=\Gal(\Kinfty/K_0)$.
\begin{lemma}\label{ghduality} Let $M$ be an perfect complex of $\Lambda=\Oh_p[[\Gamma]]$-modules.
Then there are canonical isomorphisms
\[
M^*\otimes^{\Lbb}_{\Lambda}\Oh_p\isom R\Gamma(\Gamma,M)^*
\]
where the right hand side is the (continuous) group cohomology of $\Gamma$
and $M^*=\Hom(M, \Q_p/\Z_p\otimes\Oh_p)$.
\end{lemma}
\bew
We have 
\begin{align*}
R\Hom_{\Lambda}(\Oh_p, M^*)&
=R\Hom_{\Lambda}(\Oh_p,\Hom(M,\Q_p/\Z_p\otimes\Oh_p))\\
&=R\Hom_{\Lambda}(M\otimes^{\Lbb}_{\Lambda}\Oh_p,\Q_p/\Z_p\otimes\Oh_p))
\end{align*}
which by biduality $M^{**}=M$ proves our claim.
\bewende 
\begin{cor}\label{complexes}
There are exact triangles
\begin{align*}
(\Yhinfty\otimes T_{p}E(k))\otimes^{\Lbb}_{\Lambda}\Oh_p\to& R\Gamma_{\et}
(K_0\otimes\Q_p, E[p^{\infty}](-k))^*[-1]\\
&\to R\Gamma(\Gamma, H^0_{\et}(\Kinfty\otimes
 \Q_p, E[p^{\infty}](-k)))^*[-1]
\end{align*}
and
\begin{align*}
(\Xhinfty\otimes T_{p}E(k))\otimes^{\Lbb}_{\Lambda}\Oh_p\to &R\Gamma_{\et}(\Oh_{0,S_p}, E[p^{\infty}](-k))^*[-1]\\
&\to
R\Gamma(\Gamma, H^0_{\et}(\Oh_{\infty,S_p}, [p^{\infty}](-k)))^*[-1]
\end{align*}
\end{cor}
\bew Apply lemma \ref{ghduality} to the exact triangles in 
corollary \ref{extri}.
\bewende

Now we come back to the Tate-Poitou localization
sequences over $\Oh_{0,S_p}$
\begin{align*}
R\Gamma_{\et}(\Oh_{0,S_p}, T_{p}E(k+1))\to
&R\Gamma_{\et}(K_0\otimes\Q_p, E[p^{\infty}](-k))^*[-2]\\&
\to R\Gamma_{\et}(\Oh_{0,S_p}, E[p^{\infty}](-k))^*[-2]\to
\end{align*}
and the map $e_{p}$. 
Recall that 
\[
e_{p}: \Chbarinfty\otimes T_{p}E(k)\to
H^1_{\et}(\Oh_{S_p}, T_{p}E(k+1))
\]
(see definition \ref{souleelements}) and taking in the definition of 
$e_p$ only the norm maps to $K_0$ we get a map:
\[
e_{p}: \Chbarinfty\otimes T_{p}E(k)\to
H^1_{\et}(\Oh_{0,S_p}, T_{p}E(k+1))
\]
As $H^0_{\et}(\Oh_{0,S_p}, T_{p}E(k+1))=0$ for weight reasons, 
we get a map of complexes
\[
e_{p}:(\Chbarinfty\otimes T_{p}E(k))\otimes^{\Lbb}_{\Lambda}\Oh_p\to
R\Gamma_{\et}(\Oh_{0,S_p}, T_{p}E(k+1))[1].
\]
This is compatible with the maps defined before:
\begin{lemma}\label{exactseq} The following diagram is commutative
\[
\begin{CD}
(\Chbarinfty\otimes T_{p}E(k))\otimes^{\Lbb}_{\Lambda}\Oh_p@>e_{p}>>
R\Gamma_{\et}(\Oh_{0,S_p}, T_{p}E(k+1))[1]\\
@VVV@VVV\\
(\Yhinfty\otimes T_{p}E(k))\otimes^{\Lbb}_{\Lambda}\Oh_p@>>>
R\Gamma_{\et}(K_0\otimes\Q_p,E[p^{\infty}](-k))^*[-1]\\
@VV\alpha V@VVV\\
(\Xhinfty\otimes T_{p}E(k))\otimes^{\Lbb}_{\Lambda}\Oh_p@>>>
R\Gamma(\Oh_{0,S_p},
 E[p^{\infty}](-k))^*[-1]
\end{CD}
\]
Here $\alpha $ is induced by the map $\Yhinfty/\Chbarinfty\to \Xhinfty$.
\end{lemma}
\bew The commutativity of the lower
square is clear. Let us treat the upper square.
The map 
\[
(\Yhinfty\otimes T_{p}E(k))_{\Gamma}
\to H^1_{\et}(K_0\otimes\Q_p, E[p^{\infty}](-k))^*
\]
is the dual of the corestriction
\[
H^1_{\et}(K_0\otimes\Q_p, E[p^{\infty}](-k))\to H^1_{\et}(\Kinfty\otimes
\Q_p, E[p^{\infty}](-k))^{\Gamma}.
\]
By the local duality theorem the  corestriction
map is dual to  the norm map  
\[
 H^1_{\et}(\Kinfty\otimes
\Q_p,T_{p} E(k+1))_{\Gamma}\to
H^1_{\et}(K_0\otimes\Q_p,T_{p} E(k+1)).
\]
This together with the definition of $e_{p}$ proves our claim.
\bewende

\subsubsection{The comparison theorem between elliptic units and 
Galois cohomology}

The next step is to relate the determinants of 
$(\Chbarinfty\otimes T_{p}E(k))\otimes^{\Lbb}_{\Lambda}\Oh_p$
and $R\Gamma_{\et}(\Oh_{0,S_p}, T_{p}E(k+1))$ as $\Oh_p[\Delta]$-modules. 
For this we need Rubin's ``main conjecture''. As the ``main conjecture'' is
not proven for characters of $\Delta$, which are trivial on the 
decomposition group $\Delta_p$ of $p$,  we need  the following lemma:
\begin{lemma}\label{goodchi} Let $p$ be inert or ramified, 
where $p$ is a prime over which $E$ has good reduction. 
Let $\chi$ be the $\Delta$ representation on $ \Hom_{\Oh_p}(T_{p}E(k),\Oh_p)$.
Then $\Delta_p=\Delta$ and $\chi$ is non trivial on $\Delta_p$.
\end{lemma}
\bew
Let $p$ be an inert or ramified prime
and $\chi'$ be the $\Delta$ representation on $T_{p}E$. Then
$\chi'$ is irreducible (\cite{Ru1} 11.5.) and because $E$ has good
reduction the prime above $p$ is totally ramified in $K_0$ (\cite{Ru2} 3.6.) 
and $\Delta\isom (\Oh/\pf)^*$. Now $\chi'$ is two dimensional
and $\chi$ is simply a twist of $\chi'$ by a power of $\det\chi'$. Thus 
$\chi$ acts non trivially on $ \Hom_{\Oh_p}(T_{p}E(k),\Oh_p)$.
\bewende

\begin{cor}\label{uinftyisom} Let $\chi$ be the $\Delta$ representation on 
$ \Hom_{\Oh_p}(T_{p}E(k),\Oh_p)$ and $p\nmid \Norm\ff$ be a prime. Then
\[
\Uhinfty^{\chi}\isom\Yhinfty^{\chi}.
\]
\end{cor}
\bew If $p$ is split this follows immediately from lemma \ref{UYrel}
and if $p$ is inert or prime in $K$ this follows from the same lemma
and the above result 
because the $\chi$-eigenspace of $\Oh_p[\Delta/\Delta_p]$ is zero.
\bewende
We can now formulate the main theorem of this section.

\begin{thm}\label{cohdescr}
Let $\chi$ be the $\Delta $-representation
on $\Hom_{\Oh_p}(T_{p}E(k),\Oh_p)$ and assume that  $p\nmid\Norm\ff$. 
Then the map $e_p$ induces
an isomorphism of $\Oh_p$-modules
\[
{\det}_{\Oh_p}\left( (\Chbarinfty^{\chi}\otimes_{\Oh_p}  T_{p}E(k))\otimes^{\Lbb}_{\Lambda}\Oh_p\right)\isom {\det}_{\Oh_p}\left( R\Gamma(\Oh_{S}, T_{p}E(k+1))\right)^{-1}.
\]
\end{thm}
The rest of this section is concerned with the proof of this
theorem. Let us first show:

\begin{prop}\label{vanishing} Let $\chi$ and $p$ be as in theorem
\ref{cohdescr}, then 
\begin{align*}
{\det}_{\Oh_p}(R\Gamma(\Gh, H^0_{\et}(\Kinfty\otimes\Q_p,E[p^{\infty}](-k)))&\isom \Oh_p\\
{\det}_{\Oh_p}(R\Gamma(\Gh, H^0_{\et}(\Oh_{\infty,S_p},E[p^{\infty}](-k)))&\isom \Oh_p
\end{align*}
\end{prop}
\bew The action of $\Gh$ on $T_p(k)\isom \Oh_p$ is via a character
$\Gh\to \Oh_p^*$. This gives a surjection $\Oh_p[[\Gamma]]\to T_pE(k)$.
 As $\Gamma\isom\Z_p^2$ the
kernel of this surjection is an ideal with height $2$ and hence
\[
{\det}_{\Oh_p[[\Gh]]}(T_{p}E(k))\isom \Oh_p[[\Gh]].
\]
This implies ${\det}_{\Oh_p}(T_{p}E(k)\otimes^{\Lbb}_{\Oh_p[[\Gh]]}\Oh_p)
\isom \Oh_p$.
Lemma \ref{ghduality} then implies the claim. 
\bewende
Recall that by corollary \ref{uinftyisom}  we have an isomorphism
\[
\Uhinfty^{\chi}\isom  \Yhinfty^{\chi}.
\]
\begin{cor}  The triangles in corollary \ref{complexes} give 
rise to isomorphisms
\begin{align*}
{\det}_{\Oh_p}((\Uhinfty^{\chi}\otimes_{\Oh_p} T_{p}E(k))\otimes^{\Lbb}_{\Lambda}\Oh_p)&\isom 
{\det}_{\Oh_p}( H^0(\Delta,R\Gamma_{\et}
(K_0\otimes  
\Q_p, E[p^{\infty}](-k))^*[-1]))\\
{\det}_{\Oh_p}((\Xhinfty^{\chi}\otimes_{\Oh_p} T_{p}E(k))\otimes^{\Lbb}_{\Lambda}\Oh_p)&
\isom {\det}_{\Oh_p}(H^0(\Delta,R\Gamma_{\et}(\Oh_{0,S_p}, E[p^{\infty}](-k))^*[-1]))
\end{align*}
\end{cor}
\bew The complexes in the triangle in \ref{complexes} are
$\Oh_p[\Delta]$-modules and we apply $R\Gamma(\Delta,\_ )$. Then
\[
R\Gamma(\Delta,\Yhinfty\otimes_{\Oh_p} T_{p}E(k))\isom 
\Yhinfty^{\chi}\otimes_{\Oh_p} T_{p}E(k)
\]
by definition of $\chi$. The same holds for $\Xhinfty\otimes_{\Oh_p} T_{p}E(k)$.
The result follows with proposition \ref{vanishing}.
\bewende

\begin{cor}
There is an isomorphism of determinants
\begin{eqnarray*}\lefteqn{
{\det}_{\Oh_p}\left(H^0(\Delta, R\Gamma(\Oh_{0,S_p}, T_{p}E(k+1)))\right)^{-1}
\isom }\hspace{1.5cm}\\
&&\isom {\det}_{\Oh_p}
\left((\Uhinfty^{\chi}\otimes_{\Oh_p} T_{p}E(k))\otimes^{\Lbb}_{\Lambda}\Oh_p\right)
{\det}_{\Oh_p}
\left(\Xhinfty^{\chi}\otimes_{\Oh_p} T_{p}E(k)\otimes^{\Lbb}_{\Lambda}\Oh_p\right)^{-1}
\end{eqnarray*}
\end{cor}
\bew
Apply $R\Gamma(\Delta,\_ )$ to the triangle
\begin{align*}
 R\Gamma(\Oh_{0,S_p}, T_{p}E(k+1))\to& R\Gamma_{\et}(K_0\otimes
 \Q_p, E[p^{\infty}](-k))^*[-2]\to\\
&\to 
R\Gamma(\Oh_{0,S_p},
 E[p^{\infty}](-k))^*[-2]
\end{align*}
and use the above corollary.
\bewende

Finally, we need to investigate the relation of the cohomology of $\Oh_{0,S_p}$
and $\Oh_{0,S}$, which is the integral closure of $\Oh_S$ in $K_0$.
\begin{lemma}Let $p$ and $\chi$ be as in the theorem \ref{cohdescr}.
The restriction map of cohomology of $\Oh_{0,S_p}$ to
$\Oh_{0,S}$ induces an 
equality of determinants 
\begin{align*}
{\det}_{\Oh_p}\left(H^0(\Delta,R\Gamma_{\et}(\Oh_{0,S_p}, T_{p}E(k+1)))\right)&\isom 
{\det}_{\Oh_p}\left(H^0(\Delta,R\Gamma_{\et}(\Oh_{0,S}, T_{p}E(k+1)))\right)\\
&\isom {\det}_{\Oh_p}\left(R\Gamma_{\et}(\Oh_{S}, T_{p}E(k+1)))\right).
\end{align*}
\end{lemma}
\bew 
There is an exact triangle 
\[
 R\Gamma(\Oh_{0,S_p}, T_{p}E(k+1))\to 
 R\Gamma(\Oh_{0,S}, T_{p}E(k+1))\to \bigoplus_{v\in S\ohne S_p}
R\Gamma_{\kappa(v)}(\Oh_v,T_{p}E(k+1))[1]
\]
where $\Oh_v$ is the local ring at $v$.
As $T_{p}E(k+1)$ is unramified at the places $v$ in $K_0$, which
are in $S\ohne S_p$ we have
by purity
\[
R\Gamma_{\kappa(v)}(\Oh_v,T_{p}E(k+1))\isom R\Gamma({\kappa(v)},T_{p}E(k)).
\]
Let us prove that 
\[
H^0(\Delta,\bigoplus_{v\in S\ohne S_p}R\Gamma({\kappa(v)},T_{p}E(k)))=0.
\]
For this note that $H^1({\kappa(v)},T_{p}E(k))\isom T_pE(k)_{\Gal(\bar{\kappa(v)}/{\kappa(v)})}$ are the coinvariants and that $H^0=0$. 
Fix a prime $v\in S\ohne S_p$ of $K$ dividing $\ff$, 
then the primes $v_0|v$ of $K_0$ 
are permuted by $\Delta$. Fix $v_0$ dividing $v$ and let $\Delta_{v_0}$ be
the stabilizer of $v_0$. It suffices to prove that
$\Delta_{v_0}$ acts non trivially on $T_pE(k)_{\Gal(\bar{\kappa(v)}/{\kappa(v)})}$. Let $I_{v_0}\subset\Delta_{v_0}$ be the inertia group of $v_0$.
This group is non trivial because $K_0/K$ is ramified above $v$ by lemma
\ref{redlemma}
and it acts non trivially on $T_{p}E(k)$ 
because $v_0|\ff$ and by the Neron-Ogg-Shavarevich criterium. 
This proves our claim.
\bewende

Now we can prove the theorem.
\bew(of theorem \ref{cohdescr}) Let $\chi$ and $p$  be as in the 
theorem.
By Rubin's ``main conjecture'' \ref{mainconjcor} we have
\begin{align*}
{\det}_{\Oh_p}\left((\Uhinfty^{\chi}/\Chbarinfty^{\chi})\otimes T_{p}E(k)
\otimes^{\Lbb}_{\Lambda}\Oh_p\right)&\isom
{\det}_{\Lambda}\left((\Uhinfty^{\chi}/\Chbarinfty^{\chi})\otimes T_{p}E(k)
\right)\otimes_{\Lambda}\Oh_p\\
&\isom {\det}_{\Lambda}\left(\Xhinfty^{\chi}\otimes T_{p}E(k)\right)\otimes_{\Lambda}\Oh_p\\
&\isom
{\det}_{\Oh_p}\left(\Xhinfty^{\chi}\otimes T_{p}E(k)\otimes^{\Lbb}_{\Lambda}\Oh_p\right).
\end{align*}
On the other hand, 
\[
{\det}_{\Lambda^{\chi}}(\Uhinfty^{\chi}/\Chbarinfty^{\chi})\isom
{\det}_{\Lambda^{\chi}}(\Uhinfty^{\chi})\otimes {\det}_{\Lambda^{\chi}}
(\Chbarinfty^{\chi})^{-1}.
\]
This together with the above corollaries gives the result.
\bewende

\section{The elliptic polylogarithm sheaf}
We start afresh with the aim of computing the specialization of
the elliptic polylog. For this we have to recall the definition
of the polylogarithm sheaf and give a geometric interpretation of
it.\\

We review here mostly Beilinson and Levin \cite{Be-Le}.
Everything that follows will be in the general setting of 
an elliptic curve over any base $S$. Because of this we 
start with fixing the notations. Then we review the unipotent
elliptic polylog of Beilinson and Levin. For our geometrical
construction we need a different description of this polylogarithm
sheaf in terms of the fundamental group of the elliptic curve.
This description is in section \ref{geompol}. 
The comparison of these two approaches will be carried
out  in section \ref{comparisonpol}.
Finally we consider the specialization of the polylogarithm sheaf
at torsion points. This gives the $l$-adic Eisenstein classes.

\subsection{Notations and conventions}
Let $S$ be a scheme, and $l$ be a prime number invertible on $S$.
We fix a base ring $\Lambda:=\Z/l^r\Z,\Z_l$ or $\Q_l$. In this 
section we introduce some notations for elliptic curves over $S$
and for pro-$\Lambda$-sheaves.

\subsubsection{Elliptic curves and coverings}
\begin{defn}\label{curvedef} An  {\em elliptic curve } is  a  smooth 
proper morphism $\pibar :E\to S$ together with a section $e:S\to E$,
such that the geometric fibers $E_{\sbar}$ of 
$\pibar$ are  connected curves of genus $1$. 
\end{defn}
We introduce the following {\em notation}:
On $E$ we have the multiplication by $N$ map, which we denote by $[N]$.
We let $H_n:=\ker [l^n]$ and we denote by $E_n$ the curve $E$ over 
$S$ considered as a $H_n$-torsor over $E$. The $l^n$-multiplication
map will then be denoted by $p_n:E_n\to E$. Let $U_n:=E_n\ohne H_n$ 
and $U:=E\ohne e(S)$, so that we have a
Cartesian diagram
\[
\begin{CD}
H_n@>h_n>>E_n@<j_n<<U_n\\
@VVp_{H_n}V@VVp_nV@VVV\\
S@>e>>E@<j<< U.
\end{CD}
\]
The unit section of $E_n$ will be $e_n$, if confusion is likely. 
The map $E_m\to E_n$ for $m\ge n$, which is the multiplication 
by $l^{m-n}$, is denoted by $p_{m,n}$ or even $p$. Let $\bar{\pi}_n:E_n\to S$
and $\pi_n:U_n\to S$ be the structure maps.

\subsubsection{Pro-sheaves}\label{prosheaves}

The polylogarithm is an extension of pro-sheaves and we will work
in the category of pro-sheaves. For convenience
of the reader we recall the definition and the main properties of
pro-objects in the case we need. 

Let $\Ah$ be an abelian category.
\begin{defn}
The category $\pro-\Ah$ of {\em pro-objects} 
is the category whose objects are projective systems
\[
A:I^{\op}\to \Ah
\]
denoted by $(A_i)_{i\in I}$, 
where $I$ is some small filtered index category.
The  morphisms are
\[
\Hom_{\pro-\Ah}((A_i),(B_j)):=\prolim_j\indlim_i\Hom_{\Ah}(A_i,B_j).
\]
\end{defn}
The category $\pro-\Ah$ is again abelian (see \cite{Ar-Ma} A 4.5).
We call an object $(A_i)_I\in \pro-\Ah$ {\em Mittag-Leffler zero} if
for every $i\in I$ there is an $i\to j$ such that $A_j\to A_i$ is
the zero map.
An element is zero in $\pro-\Ah$ if and only if it is Mittag-Leffler
zero (see \cite{Ar-Ma} A 3.5). A functor $F:\Ah\to \Bh$ is extended
to the pro-categories in the obvious way $F((\Fh_i)_i):=(F(\Fh_i))_i$.

Let us specialize to the category $\Sh_{\et}(X)$ of \'etale sheaves
on a scheme $X$. We denote by $\pro-\Sh_{\et}(X)$ the associated category of
pro-sheaves as defined above. 
Pro-sheaves will usually be written as $(\Fh_i)_i$ the
transition maps understood. For two pro-sheaves $(\Fh_i)_i$ and
$(\Gh_i)_i$ on a scheme $X$ define $\Ext^j_X((\Fh_i)_i,(\Gh_i)_i)$
to be the group of $j$-th Yoneda extensions of $(\Fh_i)_i$ by
$(\Gh_i)_i$ in $\pro-\Sh_{\et}(X)$.

\subsection{Review of the elliptic  polylogarithm}\label{seclogarithm}
We first recall the definition of the 
elliptic logarithm sheaf from \cite{Be-Le}.
Then we define the elliptic polylogarithmic sheaf.

\subsubsection{The unipotent logarithm sheaf}
Recall that  $\Lambda$ is either $\Z/l^r\Z$, $\Z_l$ or $\Q_l$.

\begin{defn}
A lisse $\Lambda$-sheaf $\Fh$ on $E$ is unipotent of length $n$, if 
it admits a filtration $\Fh=\Fh^0\supset\Fh^1\supset\ldots\supset\Fh^n\supset 0$ 
such that $\gr^i\Fh=\pibar^*\Gh^i $ for some lisse
$\Lambda$-sheaf $\Gh^i $ on $S$.
\end{defn}
 Let $\shH_{\Lambda}:=\underline{\Hom}_S(R^1\pibar_*\Lambda,\Lambda)$, then the boundary map
for the exact sequence $0\to \gr^{i+1}\Fh\to \Fh^i/\Fh^{i+2}\to \gr^i\Fh\to 0$
induces by duality a map
\[
\shH_{\Lambda}\otimes\gr^i\Fh\to\gr^{i+1}\Fh.
\]
 This gives an
action of the ring $S^{\le n}:=\bigoplus_{k=0}^n\Sym^k\shH_{\Lambda}$ on
$\gr^{\bullet}\Fh$. Beilinson and Levin prove:
\begin{thm}(\cite{Be-Le} 1.2.6.)\label{univlogprop}
There is a  $k$-unipotent sheaf $\Log^{(k)}$
together with a section $1^{(k)}:\Lambda\to e^*\Log^{(k)} $ of the fibre
at the unit section $e$ of $E$, which is unique up to isomorphism,
 such that for every $k$-unipotent 
sheaf $\Fh$ the map
\begin{align*}
\pibar_*\underline{\Hom}_E(\Log^{(n)},\Fh)&\to e^*\Fh\\
f&\mapsto f\verk 1^{(n)} 
\end{align*}
is an isomorphism. 
\end{thm}
Recall also from \cite{Be-Le} that this is equivalent
to the fact that the map $\nu:S^{\le k}\to \pibar_*\gr^{\bullet}\Log^{(k)}$
that sends $1$ to $1^{(n)}$ is an isomorphism. 
\begin{defn}
The canonical maps $\Log^{(k+1)}\to \Log^{(k)}$ that map $1^{(k+1)}$ to 
$1^{(k)}$ make 
\[
\Log:=(\Log^{(k)})_k
\]
 a pro-sheaf, which is called the 
{\em logarithm sheaf}. 
 If it is necessary to indicate
$\Lambda$ we write $\Log^{(k)}_{\Lambda}$ and $\Log_{\Lambda}$.
\end{defn}

Denote by $\Rh^{(k)}:=e^*\Log^{(k)}$ the fibre of $\Log^{(k)}$. This is
a ring with identity given by $1^{(k)}$. Moreover $\Rh:=e^*\Log$ has
a Hopf algebra structure. Then 
$\pibar^*\Rh^{(k)}$ acts on $\Log^{(k)}$ and for every section $t:S\to E$
the sheaf $t^*\Log^{(k)}$ is a free module of rank $1$ over $\Rh^{(k)}$.
The action of $\pibar^*\Rh^{(k)}$ on $\Log^{(k)}$ induces via the isomorphism 
$\pibar_*\underline{\Hom}_E(\Log^{(k)},\Fh)\xrightarrow{\sim} e^*\Fh$ an action of
$\Rh^{(k)}$ on $e^*\Fh$. In fact we have:
\begin{prop}(\cite{Be-Le}1.2.10 v))\label{catequiv} The map $\Fh\mapsto e^*\Fh$ 
 is an equivalence of the category of $k$-unipotent sheaves on $E$
with the category of lisse $\Rh^{(k)}$-modules on $S$.
\end{prop}
We just remark that the inverse functor
is $\Mh\mapsto \pibar^*\Mh\otimes_{\pibar^*\Rh^{(k)}}\Log^{(k)}$.

\subsubsection{Higher direct images of the logarithm sheaf}
Denote by $\Ih^{(k)}$ the augmentation ideal of the ring $\Rh^{(k)}$. 
The pro-sheaves $(\Rh^{(k)})_k$ and $(\Ih^{(k)})_k$ are
denoted by $\Rh$ and $\Ih$ respectively.
The important fact for the definition of the polylogarithm is the 
computation of the higher direct images of $\Log^{(k)}$.
\begin{prop}(see \cite{Be-Le} 1.2.7) 
The higher direct images of $\Log^{(k)}_{\Lambda}$ are
\[
R^i\pibar_*\Log^{(k)}_{\Lambda}=\left\{\begin{array}{lll}
\Sym^k\shH_{\Lambda}&\mbox{ if }& i=0\\
\Sym^{k+1}\shH_{\Lambda}(-1)&\mbox{ if }& i=1\\
\Lambda(-1)&\mbox{ if }& i=2
\end{array}\right.
\]
The transition maps $R^i\pibar_*\Log^{(k+1)}\to R^i\pibar_*\Log^{(k)}$
are zero for $i=0,1$ and the identity for $i=2$. In particular 
$ R^i\pibar_*\Log=0$ for $i=0,1$ and $ R^2\pibar_*\Log=\Lambda(-1)$.
\end{prop}
For all the properties of the logarithm sheaf
we refer to section 1.2. in \cite{Be-Le}. 

\begin{rem} Note that in the
case of $\Lambda=\Q_l$ we have an isomorphism $\Log^{(k)}\isom \Sym^k\Log^{(1)}$
which sends $1^{(k)}$ to $1^{(1)k}/k!$. This  approach to the logarithm sheaf
is used in \cite{Hu-Ki}.
\end{rem}

Recall that $U:=E\ohne e$  is the complement of the unit section and 
$\pi:U\to S$ its structure map.
\begin{prop} \label{logimages}
The pro-sheaves $(R^i\pi_*\Log^{(k)})_k$ are Mittag-Leffler
zero for $i\neq 1$ and the canonical map
\[
R^1\pi_*\Log^{(k)}\to e^*\Log^{(k)}(-1)=\Rh^{(k)}(-1)
\]
induces an isomorphism of pro-sheaves $(R^1\pi_*\Log^{(k)}(1))_k\isom(\Ih^{(k)})_k$.
\end{prop}
\bew
Consider the localization sequence
\[
\to R^i\pibar_*\Log^{(k)}\to R^i\pi_*\Log^{(k)}\to R^{i+1}e^!\Log^{(k)}\to
\]
and the purity isomorphism $R^{2}e^!\Log^{(k)}=e^*\Log^{(k)}(-1)$. Moreover 
$R^{i}e^!\Log^{(k)}=0$ for $i\neq 2$. This together with the above
values of $R^i\pibar_*\Log$ gives the desired result.
\bewende

\subsubsection{The polylogarithm sheaf}
We are now going to define the elliptic polylogarithm. 
We will not use the usual approach using an identification of an $\Ext$ 
with a $\Hom$-group but an other direct construction due to 
Beilinson and Levin \cite{Be-Le} 1.3.6. This has the advantage
of giving directly a pro-sheaf and not only an extension class.
Moreover this sheaf can be easily compared to the geometric construction
we give later. 

For every sheaf $\Fh$ on $E$ we denote by $\Fh_U$ its
restriction to $U$. 
Let $\Fh$ be a lisse $\Lambda$-sheaf on $E$ and consider the open immersion
\[
j:U\times_SU\ohne \Delta\hookrightarrow U\times_SU
\]
where $\Delta $ is the diagonal. Define a lisse $\Lambda$-sheaf
$H_e(\Fh)$ on $U$ as follows:
\begin{defn}\label{Hdefn}Define a functor from lisse $\Lambda$-sheaves on 
$E$ to lisse $\Lambda$-sheaves on $U$ by
\[
H_e(\Fh):=R^1\pr_{1*}j_!\pr^*_2\Fh_U,
\]
where $\pr_2$ is the projection of $U\times_SU\ohne \Delta$ to the second
factor.
\end{defn}
The exact sequence
\[
0\to j_!\pr^*_2\Fh_U\to \pr^*_2\Fh_U\to \Delta_*\Fh_U\to 0
\]
induces an exact sequence
\begin{equation}\label{polseq}
0\to \pi^*\pibar_*\Fh\to \Fh_U\xrightarrow{\alpha} 
H_e(\Fh)\to \pi^*R^1\pi_*\Fh_U\to 0.
\end{equation}
Obviously this sequence is functorial in $\Fh$ so that we have the same 
sequence for pro-sheaves $(\Fh_k)_k$.
\begin{lemma} The sequence (\ref{polseq}) induces an exact sequence of
pro-sheaves
\[
0\to \Log(1)_U\to H_e(\Log(1))\to \pi^*\Ih \to 0.
\]
\end{lemma}
\bew By proposition \ref{logimages} we have $\pibar_*\Log(1)=0$ and
$R^1\pi_*\Log(1)\isom \Ih$. This implies the claim.
\bewende

Thus $H_e(\Log(1))$ gives a class in $\Ext^1_U( \pi^*\Ih,\Log(1)_U)$
which was defined in section
\ref{prosheaves} as the group of Yoneda extensions in the
abelian category $\pro-\Sh(U)$. We have 
$ H_e(\Log(1)\otimes \pi^*\Rh)\isom  H_e(\Log(1))\otimes \pi^*\Rh$ 
so that the action of $\pi^*\Rh$ on $\Log$ gives a $\pi^*\Rh$-module
structure on $H_e(\Log(1))$. In particular $H_e(\Log(1))$ is
a class in $\Ext^1_{U,\pi^*\Rh}( \pi^*\Ih,\Log(1)_U)$, i.e. a Yoneda
extension of $\pi^*\Rh$-modules.
\begin{defn}\label{poldefn} The pro-sheaf
\[
\pol:=H_e(\Log(1))
\]
is the {\em elliptic polylogarithm sheaf}. If we need to indicate the
dependence on $\Lambda$ we write $\pol_{\Lambda}$. We also define
$\pol^{(k)}:=H_e(\Log^{(k)}(1))$.
\end{defn}

\subsection{A geometric approach to the elliptic polylog sheaf}\label{geompol}
We now present a different construction of the logarithm sheaf in
the case $\Lambda=\Z/l^r\Z$.
This makes explicit the remark in \cite{Be-Le} 1.2.5.

\subsubsection{The geometric logarithm sheaf} 
Recall that $p_n=[l^n]:E_n\to E$ and consider the sheaves 
\[
\Logeom_n:={p_n}_*\Lambda
\]
on $E$. For $m\ge n$ we have the trace map ${p_{m}}_*\Lambda\to {p_n}_*\Lambda$
and we define:
\begin{defn}
The {\em geometric logarithm  sheaf} is the pro-sheaf
\[
\Logeom:=(\Logeom_n)_n
\]
where the transition maps are the 
above trace maps. Let 
\[
\Rhgeom:=(\Rhgeom_n)_n:= (e^*\Logeom_n)_n
\]
be the pro-sheaf defined by the
pull-back of $\Logeom_n$ along the unit section $e$.
Let $\Ihgeom:=(\Ihgeom_n)_n:=\ker(\Rhgeom\to \Lambda)$ be the
augmentation ideal of $\Rhgeom$.  
\end{defn}
Note that the existence of the section $e_n$ of $H_n=p_n^{-1}(e)$
implies that there is a map $1_n:\Lambda\to e^*\Logeom_n=\Rhgeom_n$.
The action of $H_n$ on $E_n$ over $E$ gives an action of
$H_n$ on $\Logeom_n$, hence an action of $\pibar^*\Rhgeom_n$ on $\Logeom_n$.

The sheaf $\Logeom_n$ has the following important property. Recall
that $\Lambda=\Z/l^r\Z$.
\begin{prop}\label{geomunivprop} For every lisse $\Lambda$-sheaf $\Fh$
the map 
\begin{align*}
\indlim_n\pibar_*\underline{\Hom}_E(\Logeom_n,\Fh)&\to e^*\Fh\\
f&\mapsto f\verk 1_n
\end{align*}
is an isomorphism. 
\end{prop}
\bew
The map $p_n$ is finite \'etale, so that
\[
\underline{\Hom}_E(p_{n*}\Lambda,\Fh)=\underline{\Hom}_{E_n}(\Lambda,p^*_n\Fh)
=p^*_n\Fh.
\]
As $\Fh$ is a lisse $\Lambda$-sheaf, there is an $n$ such that 
$p^*_n\Fh$ comes from $S$, i.e. $p^*_n\Fh\isom \pibar^*_ne^*_np^*_n\Fh$.
Thus
\[
\indlim_n\pibar_*\underline{\Hom}_E(\Logeom_n,\Fh)=\indlim_n\pibar_{n*}\pibar^*_ne^*_np^*_n\Fh=
\indlim_ne^*_np^*_n\Fh=e^*\Fh,
\]
which proves our claim.
\bewende

Let $\Fh$ be a lisse $\Lambda$-sheaf, then the action of $\Rhgeom_n$ on
$\Logeom_n$ induces via the above proposition
an action of $\Rhgeom_n$ on $e^*\Fh$ for
some $n$.
\begin{cor} \label{categoryequiv}
The functor $\Fh\mapsto e^*\Fh$ induces an equivalence of
the category of lisse $\Lambda$-modules on $E$ and lisse 
$\Lambda$-modules on $S$ with a continuous action of the 
pro-sheaf $(\Rhgeom_n)_n$ (i.e. the action factors through $\Rhgeom_m$
for some $m$).
\end{cor}
\bew The inverse functor is given by 
\[
\Mh\mapsto \pibar^*\Mh\otimes_{\pibar^*\Rhgeom_m}\Logeom_m
\]
if the action of $(\Rhgeom_n)_n$ factors through $\Rhgeom_m$.
\bewende

\subsubsection{The higher direct images of the geometric logarithm sheaf}
As for $\Log $ we  can compute the higher direct images of $\Logeom$:
\begin{lemma}
The pro-sheaf 
\[
(R^i\pibar_*\Logeom_n)_n
\] 
is Mittag-Leffler zero for $i\neq 2$ and 
\[
(R^{2}\pibar_*\Logeom_n)_n\isom \Lambda(-1).
\] 
\end{lemma}
\bew We have to compute the transition maps in 
\[
(R^i\pibar_*{p_n}_*\Lambda)_n=(R^i\pibar_{n*}\Lambda)_n\isom
(R^i\pibar_{*}\Lambda)_n
\]
where now the transition maps $R^i\pibar_{*}\Lambda\to R^i\pibar_{*}\Lambda$ 
are given by multiplication with $(l^{m-n})^{2-i}$. This map is zero for 
$m\ge n+r$. 
\bewende
\begin{cor} 
The pro-sheaves $(R^i\pi_*\Logeom_n)_n$ are Mittag-Leffler
zero for $i\neq 1$ and the canonical map
\[
R^1\pi_*\Logeom_n\to e^*\Logeom_n(-1)=\Rhgeom_n
\]
induces an isomorphism of pro-sheaves $(R^1\pi_*\Logeom_n)_n\isom(\Ihgeom_n)_n$.
\end{cor}
\bew See the proof of \ref{logimages}.
\bewende

\subsubsection{The geometric polylogarithm sheaf}
The geometric polylog sheaf can now be defined in the same way as in
\ref{poldefn}. Recall the functor $H_e$ from definition \ref{Hdefn}.
\begin{defn}\label{polgeomdefn} 
The {\em geometric elliptic polylog sheaf} is the pro-sheaf
\[
\polgeom:=H_e(\Logeom(1)).
\]
To indicate the dependence on $\Lambda$ we write $\polgeom_{\Lambda}$ and
we define 
\[
\polgeom_n:=H_e(\Logeom_n(1)).
\]
\end{defn}
As before the pro-sheaf $\polgeom$ is a sheaf of $\pi^*\Rhgeom$-modules
and defines a Yoneda extension class
in 
\[
\Ext^1_{U,\pi^*\Rhgeom}(\pi^*\Ihgeom,\Logeom(1)).
\]

\subsection{The comparison of $\pol$ and $\polgeom$}\label{comparisonpol}
In this section we compare $\pol$ and $\polgeom$. 
Recall that $\polgeom$ is only defined for $\Lambda=\Z/l^r\Z$. 

\subsubsection{Comparison of the logarithm sheaves}
We  first  reformulate the property of being a continuous $\Rhgeom$-module
on $S$. 
 Choose a geometric point $\sbar\in S$ and define a ring
\begin{equation}\label{groupringdefn}
\Lambda[[\shH_{\Z_l,\sbar}]]:=\prolim_n(\Rhgeom_n)_{\sbar}=
\prolim_n\Lambda[H_{n,\sbar}].
\end{equation}
Then a lisse $\Lambda$-sheaf $\Mh$ on $S$ with an action of $\Rhgeom_n$ 
for some $n$ is the same as a finite $\Lambda$-module with a  continuous
action of $\Lambda[[\shH_{\Z_l,\sbar}]]$. We call these continuous 
$\Lambda[[\shH_{\Z_l,\sbar}]]$-modules.
We define two ideals in $\Lambda[[\shH_{\Z_l,\sbar}]]$:
\begin{align}
\Jf_n&:=\ker(\Lambda[[\shH_{\Z_l,\sbar}]]\to \Rhgeom_n)\nonumber\\
\af&:=\ker(\Lambda[[\shH_{\Z_l,\sbar}]]\to \Lambda).
\end{align}

The universal property \ref{geomunivprop} of $\Logeom$ and the section
$1^{(k)}$ of $e^*\Log^{(k)}_{\Lambda}$ implies:

\begin{lemma} Let $\Lambda=\Z/l^r\Z$, then there is a unique map
\[
\rho^{(k)}\in\indlim_n\pibar_*\underline{\Hom}_E(\Logeom_n,\Log^{(k)}_{\Lambda})
\]
corresponding to $1^{(k)}:\Lambda\to e^*\Log^{(k)}_{\Lambda}$.
\end{lemma}
In fact we want to show that $\rho^{(k)}$ induces an isomorphism
of pro-sheaves $\rho:(\Logeom_n)_n\to (\Log^{(k)}_{\Lambda})_k$.

\begin{lemma}\label{groupringlemma} Let $\Lambda=\Z/l^r\Z$ and 
let $\af^k$ be the $k$-th power of $\af$, then 
$\af^k/\af^{k+1}\isom \Sym^k(\shH_{\Z_l}\otimes\Lambda)$ and there is 
a canonical surjection
\[
\rho^{(k)}_{r+k-2}:\Rhgeom_{r+k-2,\sbar}=\Lambda[[\shH_{\Z_l,\sbar}]]/\Jf_{r+k-2}\to \Lambda[[\shH_{\Z_l,\sbar}]]/\af^k
\]
for all $k\ge 1$, which in the limit gives an isomorphism of 
$\Lambda[[\shH_{\Z_l,\sbar}]]$ to itself.
\end{lemma}
\bew Let $\sbar\in S$ be a geometric point. Then we choose two topological
generators $\gamma, \eta$ of $\shH_{\Z_l,\sbar}$. It is a standard fact in Iwasawa theory
that we have a ring isomorphism $\Lambda[[\shH_{\Z_l,\sbar}]]\isom\Lambda[[X,Y]]$,
mapping $\gamma -1\mapsto X$ and $\eta-1\mapsto Y$, with the
power series ring in two variables. This implies that the
augmentation ideal $\af$ is mapped to the ideal $(X,Y)$. The
claim that $\af^k/\af^{k+1}\isom \Sym^k(\shH_{\Z_l}\otimes\Lambda)$
follows immediately.
Now $\Jf_k$ corresponds under this isomorphism to the ideal 
$((X+1)^{l^k}-1,(Y+1)^{l^k}-1)$. By induction one sees that this ideal is
contained in $(l,X,Y)^{k+1}$. As $l^r=0$ in $\Lambda$ we have
$(l,X,Y)^{r+k}\subset (X,Y)^{k+1}$. This gives the map as
indicated, which in the limit is clearly an isomorphism. 
\bewende
Denote by $G^{(k)}$ the sheaf on $E$ defined by the $\Rhgeom_{r+k-2}$-module
$\Lambda[[\shH_{\Z_l,\sbar}]]/\af^k$ via the equivalence of categories in
\ref{categoryequiv}. This is a unipotent sheaf with a section $1$ of $e^*G^{(k)}$
given by the identity in $\Lambda[[\shH_{\Z_l,\sbar}]]/\af^k$.
Thus, by \ref{univlogprop} there is a unique map 
$\Log^{(k)}\to G^{(k)}$ mapping $1^{(k)}$ to $1$.

\begin{prop} \label{logisom}Let $\Lambda=\Z/l^r\Z$.
The map $\Log^{(k)}\to G^{(k)}$ is an isomorphism. In particular 
the map of pro-sheaves $\rho:\Logeom_{\Lambda}\to \Log_{\Lambda}$
defined above is an isomorphism.
\end{prop}
\bew The surjection $\Logeom_{r+k-2}\to G^{(k)} $ factors through
$\Log^{(k)}\to G^{(k)}$ so that this map is also surjective. 
As these two sheaves are finite of the same cardinality, this is also
an isomorphism. 
\bewende
\subsubsection{Comparison with the $\Z_l$--version}
Let $\Lambda_r:=\Z/l^r\Z$ and consider the pro-sheaf in $r$ and
$n$
\[
(\Logeom_{\Lambda_r,n})_{r,n},
\]
where $\Logeom_{\Lambda_{r+1},n}\to \Logeom_{\Lambda_r,n}$ is
induced by the reduction map $\Lambda_{r+1}\to \Lambda_r$.
We observe:
\begin{lemma} For every $k\ge 0$, 
the projective system $(\Logeom_{\Lambda_r,r+k-2})_{r}$
is cofinal in $(\Logeom_{\Lambda_r,n})_{r,n}$.
\end{lemma}

The $\Z_l$-version of $\Log$ is reconstructed from the $\Lambda_r$--version 
of $\Logeom$ as follows:
\begin{prop}\label{logcomparison}
The map $\rho^{(k)}$ induces an morphism of pro-sheaves
\[
(\Logeom_{\Lambda_r,r+k-2})_{r}\to (\Log^{(k)}_{\Lambda_r})_r=\Log^{(k)}_{\Z_l}.
\]
The induced morphism of pro-sheaves
\[
\rho:(\Logeom_{\Lambda_r,r})_{r}\isom(\Logeom_{\Lambda_r,r+k-2})_{r,k}\to (\Log^{(k)}_{\Z_l})_k=\Log_{\Z_l}
\]
is an isomorphism.
\end{prop}
\bew From lemma \ref{groupringlemma} and the universal property of 
$\Logeom$ we get the morphism of $(\Logeom_{\Lambda_r,r+k-2})_{r}$
to $\Log^{(k)}_{\Z_l}$. To check that this induces an isomorphism 
if take pro-sheaves in the index $k$, it is enough to 
check this after pull-back with $e^*$. We get the map
\[
\prolim_k\prolim_r\Lambda_r[H_{r+k-2,\sbar}]\to \prolim_k\Z_l[[\shH_{\Z_l,\sbar}]]/\af^k,
\]
which is an isomorphism.
\bewende
\subsubsection{Comparison of polylogarithm sheaves}
Using the functor $H_e$ from \ref{Hdefn} we can translate the comparison results
for the logarithm sheaves to the polylog.
\begin{prop}
The isomorphism  $\rho$ form \ref{logcomparison} induces an isomorphism of pro-sheaves
\[
H_e(\rho):(\polgeom_{\Lambda_r,r})_r=
(H_e(\Logeom_{\Lambda_r,r}(1)))_r\xrightarrow{\isom}\pol_{\Z_l}.
\]
\end{prop}
\bew Clear from the definition.
\bewende

\subsection{Specialization of the elliptic polylogarithm sheaf, $l$-adic Eisenstein classes}\label{specsec}

The specialization along torsion sections
of the elliptic polylog gives interesting cohomology classes. 
These are the $l$-adic
Eisenstein classes investigated in \cite{Hu-Ki1}, \cite{Hu-Ki}.
We recall their construction.

\subsubsection{Invariance of the logarithm sheaf under translation by torsion sections}
Let $N\in\Z$ be invertible on $S$ and $[N]:E\to E$ the $N$-multiplication.
The universal property \ref{univlogprop} gives us  canonical maps sending $1^{(k)}$
to $[N]^*1^{(k)}$
\[
\Log_{\Lambda}^{(k)}\to [N]^*\Log_{\Lambda}^{(k)}
\]
for $\Lambda=\Z/l^r\Z, \Z_l$ or  $\Q_l$.
Thus for every $N$-torsion point $t:S\to E$ we get a map of pro-sheaves
\[
\pr^N_t:t^*\Log_{\Lambda}\to e^*\Log_{\Lambda}.
\]
Similarly, we have a map 
\[
\pr^N_t:\Logeom_n\to [N]^*\Logeom_n
\]
for the geometric logarithm sheaf.
\begin{lemma}\label{translation}
If $l\nmid N$, then the map $\pr^N_t$ is an isomorphism.
\end{lemma}
\bew Let us first treat the unipotent
case. It suffices to show that $\Log_{\Lambda}\to [N]^*\Log_{\Lambda}$
is an isomorphism. Using the equivalence of
categories \ref{catequiv}, this can be tested after pull-back with $e^*$.
The resulting map $\Rh_{\Lambda}\to\Rh_{\Lambda}$ is on the 
the $k$-th graded piece $\Sym^k\shH_{\Lambda}$ induced by the 
map $\shH_{\Lambda}\xrightarrow{[N]^*}\shH_{\Lambda}$, which is just
the $N$-multiplication and thus an isomorphism.
In the geometric case we have $t^*\Logeom_n=\Lambda[p_n^{-1}(t)]$
and $e^*\Logeom_n=\Lambda[H_n]$. The map $\pr^N_t$ is by definition
induced by $t_n\mapsto [N]t_n$ for $t_n\in p_n^{-1}(t)$. This is 
obviously an isomorphism if $l\nmid N$.
\bewende
To define a morphism $\pr_t$ independent of $N$, we let
\begin{equation}\label{prdefn}
\pr_t:=\pr^N_t\verk(\pr^N_e)^{-1}:t^*\Log_{\Lambda}\to e^*\Log_{\Lambda}.
\end{equation}
We need an explicit description of the map $\pr_t$ on the 
geometric logarithm sheaf. 

We have $t^*\Logeom_n=\Lambda[p_n^{-1}(t)]$ and $e^*\Logeom_n=\Lambda[H_n]$.
Here $p_n^{-1}(t)\subset E[Nl^n]$ and if $l\nmid N$ we have
$E[Nl^n]=E[N]\oplus E[l^n]$.
\begin{lemma}The map from (\ref{prdefn})
\[
\pr_t :\Lambda[p_n^{-1}(t)]\to \Lambda[H_n]
\]
is induced by the projection
of $t_n\in p_n^{-1}(t)\subset E[N]\oplus E[l^n] $ to $H_n=E[l^n]$. 
\end{lemma}
\bew The map $\pr^N_t:\Lambda[p_n^{-1}(t)]\to \Lambda[H_n]$ maps
$t_n\in p_n^{-1}(t)$ to $[N]t_n\in H_n$. This gives the result.
\bewende
Passing to the limit we get a map of pro-sheaves
\[
\pr_t :\Z_l[[\shH_{\Z_l,t}]]\to \Z_l[[\shH_{\Z_l}]]
\]
where $\shH_{\Z_l,t}:=\prolim_n p_n^{-1}(t)$. This is 
induced by the projection $\shH_{\Z_l,t}\to \shH_{\Z_l}$.

\subsubsection{The moment map}
The pro-sheaf $e^*\Log_{\Q_l}=\Rh_{\Q_l}$ 
is a Hopf algebra (see \cite{Be-Le} 1.2.10 iv)) and the map
$\nu:\Sym^{\le k}\shH_{\Q_l}\to \gr^{\le k}\Rh^{(k)}_{\Q_l}$ is an isomorphism.
Let $\hat{\Uf}(\shH_{\Q_l})$ be the completion of the universal 
enveloping algebra of the abelian Lie algebra $\shH_{\Q_l}$. The canonical
filtration makes this a pro-sheaf  $(\Uf^k(\shH_{\Q_l}))_k$.

The structure theorem, \cite{Bour} ch. II, paragraph 1, no. 6, gives:
\begin{lemma} The map $\nu:\shH_{\Q_l}\to \Rh_{\Q_l}$ extends
to an isomorphism of Hopf algebra pro-sheaves
\[
\nu:\hat{\Uf}(\shH_{\Q_l})\isom  \Rh_{\Q_l},
\]
which on the $k$-th graded piece $\Sym^k\shH_{\Q_l}$ is multiplication
by $k!$.
\end{lemma}

Recall from \ref{logcomparison} that we have an isomorphism 
\[
\rho:(\Rhgeom_{\Lambda_r})_r\isom \Rh_{\Z_l}.
\]
\begin{defn} We define the {\em $k$-th moment map} $\mu^k$ to be the composition
\[
\mu^k:(\Rhgeom_{\Lambda_r})_r\isom \Rh_{\Z_l}\to \Rh_{\Q_l}\isom \hat{\Uf}(\shH_{\Q_l})\to \Sym^k\shH_{\Q_l},
\]
where the last map is the projection onto the $k$-th factor.
\end{defn}

We need an explicit description of the moment map. For this we write
$(\Rhgeom_{n+k-2,\Lambda_n})_n=(\Rhgeom_{\Lambda_r})_r$ and we describe 
$\mu^k$ on $\Lambda_n[H_{n+k-2}]=\Rhgeom_{n+k-2,\Lambda_n}$.
\begin{lemma}\label{momentmap}
Write $\Sym^k\shH_{\Q_l}=\Sym^k\shH_{\Z_l}\otimes\Q_l$, then
\begin{align*}
\mu^k: (\Lambda_n[H_{n+k-2}])_n&\to  (\Sym^k H_n)_n\otimes\Q_l\\
(\sum_{h\in H_{n+k-2}} n_h(h))_n&\mapsto (\sum_{h\in H_n} n_h( h^{\otimes k}))_n\otimes\frac{1}{k!}.
\end{align*}
\end{lemma}
\bew
This follows immediately from the definition of the maps $\mu^k$ and 
$\rho$.
\bewende

\subsubsection{A splitting}

\begin{lemma}\label{map} Let $\Lambda=\Z/l^r\Z$, $\Z_l$ or $\Q_l$.
Then the inclusion $\Ih_{\Lambda}\hookrightarrow\Rh_{\Lambda}$ induces
an injective map
\[
\Rh_{\Lambda}(1)=\underline{\Hom}_{S,\Rh_{\Lambda}}(\Rh_{\Lambda},\Rh_{\Lambda}(1))
\to\underline{\Hom}_{S,\Rh_{\Lambda}}(\Ih_{\Lambda},\Rh_{\Lambda}(1)),
\]
which is an isomorphism for $\Lambda=\Q_l$.
\end{lemma}
\bew
This follows from $\underline{\Hom}_{S,\Rh_{\Lambda}}({\Lambda},\Rh_{\Lambda}(1))=0$ and 
$\underline{\Ext}_{S,\Rh_{\Q_l}}^1({\Q_l},\Rh_{\Q_l}(1))=0$, which is a consequence of the Kozsul resolution 
as $\Rh_{\Q_l}$ is a sheaf of regular rings.
\bewende
\begin{cor}\label{injection}
The local to global spectral sequence for $\Ext$ gives an injection
\[
\Ext^1_{S,{\Q_l}}({\Q_l},\Rh_{\Q_l}(1))\isom\Ext^1_{S,\Rh_{\Q_l}}(\Rh_{\Q_l},\underline{\Hom}_{S,\Rh_{\Q_l}}(\Ih_{\Q_l},\Rh_{\Q_l}(1)))\xrightarrow{a}
\Ext^1_{S,\Rh_{\Q_l}}(\Ih_{\Q_l},\Rh_{\Q_l}(1)).
\]
\end{cor}
Using the isomorphism $\nu:\Rh_{\Q_l}\isom \hat{\Uf}(\shH_{\Q_l})$ 
we have an exact sequence (Kozsul resolution for the Lie algebra $\shH_{\Q_l}$)
\[
0\to \Rh_{\Q_l}(1)\to \shH_{\Q_l}\otimes_{\Q_l}\Rh_{\Q_l}\xrightarrow{b} \Ih_{\Q_l}\to 0.
\]
Here we used $\Lambda^2\shH_{\Q_l}\isom\Q_l(1)$ induced by the Weil-pairing. 
This map $b$ induces 
\[
\begin{CD} \Ext^1_{S,{\Q_l}}({\Q_l},\Rh_{\Q_l}(1))@>a >>
\Ext^1_{S,\Rh_{\Q_l}}(\Ih_{\Q_l},\Rh_{\Q_l}(1))\\
@III@VVb^* V\\
@EEE \Ext^1_{S,\Rh_{\Q_l}}(\shH_{\Q_l}\otimes_{\Q_l}\Rh_{\Q_l},\Rh_{\Q_l}(1))\\
@III@VV=V\\
@EEE\Ext^1_{S,{\Q_l}}(\shH_{\Q_l},\Rh_{\Q_l}(1)).
\end{CD}
\]
\begin{lemma}\label{splitting} The map 
\[
 \Ext^1_{S,{\Q_l}}({\Q_l},\Rh_{\Q_l}(1))\xrightarrow{b^*\verk a}
\Ext^1_{S,{\Q_l}}(\shH_{\Q_l},\Rh_{\Q_l}(1))
\]
has a canonical splitting.
\end{lemma}
\bew We have an isomorphism 
\[
\Ext^1_{S,{\Q_l}}(\shH_{\Q_l},\Rh_{\Q_l}(1))\isom
\Ext^1_{S,{\Q_l}}(\Q_l, \underline{\Hom}_{S,\Q_l}(\shH_{\Q_l},\Rh_{\Q_l}(1)))
\]
and a contraction map 
\[
\underline{\Hom}_{S,\Q_l}(\shH_{\Q_l},\Rh_{\Q_l}(1)))\to
\Rh_{\Q_l}(1)
\]
given on
\[
\underline{\Hom}_{S,\Q_l}(\shH_{\Q_l},\Sym^k\shH_{\Q_l}(1))\isom\underline{\Hom}_{S,\Q_l}(\shH_{\Q_l},\Q_l)\otimes\Sym^k\shH_{\Q_l}(1)
\]
 by
\begin{align}
\underline{\Hom}_{S,\Q_l}(\shH_{\Q_l},\Q_l)\otimes\Sym^k\shH_{\Q_l}(1)&\to
\Sym^{k-1}\shH_{\Q_l}(1)\\
f\otimes h_1\otimes\ldots\otimes h_k&\mapsto \frac{1}{k+1}\sum_{i=1}^k f(h_i)h_1\otimes
\ldots \hat{h}_i\ldots\otimes h_k.\nonumber
\end{align}
This gives the required  map  and it is straightforward to check that this is
indeed a splitting of ${b^*\verk a}$.
\bewende
\subsubsection{The specialization of the polylogarithm, $l$-adic Eisenstein classes}
We now define the specialization of the elliptic 
polylogarithm. Let $\beta\in \Z[E[N](S)\ohne e]$ be 
of the form 
\[
\beta=\sum_{t\in E[N](S)\ohne e}n_tt.
\]
We want to define an element 
\[
(\beta^*\pol_{\Q_l})^{k}\in H^1_{\et}(S, \Sym^k\shH_{\Q_l}(1)).
\]
First observe that $\pr_tt^*\pol_{\Q_l}$ is an element in
\[
\Ext^1_{S,\Rh_{\Q_l}}(\Ih_{\Q_l},t^*\Log_{\Q_l}(1))\xrightarrow{\pr_t}
\Ext^1_{S,\Rh_{\Q_l}}(\Ih_{\Q_l},e^*\Log_{\Q_l}(1))=
\Ext^1_{S,\Rh_{\Q_l}}(\Ih_{\Q_l},\Rh_{\Q_l}(1)).
\]
Define 
\[
\sigma: \Ext^1_{S,\Rh_{\Q_l}}(\Ih_{\Q_l},\Rh_{\Q_l}(1))\to
\Ext^1_{S,{\Q_l}}({\Q_l},\Rh_{\Q_l}(1))
\]
to be the composition of $b^*$ and the splitting of lemma \ref{splitting}. 
This is
a splitting of $a$. Denote by $\sigma^k$ the projection onto 
\[
\Ext^1_{S,{\Q_l}}({\Q_l},\Sym^k\shH_{\Q_l}(1))=H^1({S},\Sym^k\shH_{\Q_l}(1)). 
\]
\begin{defn}\label{polylogspecialization}
 For $\beta=\sum_{t\in E[N](S)\ohne e}n_tt\in \Z[E[N](S)\ohne e]$
define 
\[
(\beta^*\pol_{\Q_l})^{k}:=\sum_{t\in E[N](S)\ohne e}n_t(\sigma^k\pr_tt^*\pol_{\Q_l})\in H^1({S},\Sym^k\shH_{\Q_l}(1))
\]
These are the $l$-adic {\em Eisenstein classes} associated to $\beta$.
\end{defn}
\begin{rem} This  numbering disagrees with the one of \cite{Hu-Ki}, 
where we  wrote $(\beta^*\pol_{\Q_l})^{k+1}$ for this class. 
\end{rem}

Our aim is to compute not the specialization but a variant of it.
Namely let $[\af]:E\to E$ be an isogeny of degree $\Norm\af:=\deg[\af]$
prime to $l$
(the notations are chosen to fit the CM case which is our ultimate 
goal). Consider the operator $([\af]^*-\Norm\af)$ on
$H^1({S},\Sym^k\shH_{\Q_l}(1))$. Here $[\af]^* $ is the pull-back
with $[\af]$ and $\Norm\af$ is the multiplication by the number 
$\deg[\af]\in\Z$. This acts on $H^1({S},\Sym^k\shH_{\Q_l}(1))$
as $\af^k\Norm\af-\Norm\af$, which is not zero if $k\neq 0$.
Thus for  $k\neq 0$ it is the same to compute 
\[
([\af]^*-\Norm\af)(\beta^*\pol_{\Q_l})^{k}\in H^1({S},\Sym^k\shH_{\Q_l}(1))
\]
or $(\beta^*\pol_{\Q_l})^{k}$.

\section{The $l$-adic realization of the elliptic polylog}\label{seccomp}
This  part is concerned with the construction of 
 the polylog on elliptic curves
in a geometric way,
which allows to compute its specializations explicitly.
This is the technical heart of the paper and in our opinion
our main contribution to the problem of the Tamagawa number 
conjecture for elliptic curves. 
\subsection{The polylog as a one-motive}
From the definition of the polylogarithm it is quite obvious that
it is a pro-sheaf consisting of $\Z/l^r\Z$-realizations of
one-motives. In this section we make
this connection more explicit. 

\subsubsection{A reformulation}\label{reformulation}

Let $\Lambda_r=\Z/l^r\Z$ and
recall from (\ref{polseq}) that the geometric polylog sits in an exact sequence
\[
0\to \pi^*\pibar_*\Logeom_n(1)\to \Logeom_n(1)_U\to \polgeom_n \to 
\pi^*R^1\pi_*\Logeom_{n,U}(1)\to 0.
\]
By the definition of the geometric logarithm $\Logeom_n=p_{n*}\Lambda_r$
and hence $\pibar_*\Logeom_n=\Lambda_r$.  We get
\[
0\to\Lambda_r(1)\to p_{n*}\Lambda_{r,U}(1)\to \polgeom_n \to 
\pi^*R^1\pi_{n*}\Lambda_{r,U}(1)\to 0.
\]

For a lisse sheaf $\Fh$ on $U$ consider the dual
\[
(\Fh)^{\vee}:=\underline{\Hom}_U(\Fh,\Lambda).
\]
Then, using Poincar\'e duality,  we get by dualizing and twisting by $1$
an exact sequence
\[
0\to\pi^*_nR^1\pi_{n!}\Lambda_{r,U}(1)\to(\polgeom_{n})^{\vee}(1)
\to  p_{n*}\Lambda_{r,U}\to \Lambda_{r,U}\to 0.
\]
Denote by $I_{n,\Lambda_r}$ the kernel of the map 
$p_{n*}\Lambda_{r,U}\to \Lambda_{r,U}$, then $(\polgeom_{n})^{\vee}(1)$
gives a class in 
\[
\Ext^1_U(I_{n,\Lambda_r}, \pi^*_nR^1\pi_{n!}\Lambda_{r,U}(1)).
\]
We want to give a geometric interpretation of this class.
For this we will relate $(\polgeom_{n})^{\vee}(1)$
with the $l^r$-torsion points of a one-motive, which is defined 
via a generalized Picard scheme.

\subsubsection{The generalized Picard scheme}\label{picardscheme}
We will give a geometric interpretation of $(\polgeom_n)^{\vee}(1)$.
For this we need the Picard scheme of line bundles on 
$E_n$ trivialized along $H_n$ (cf. the article by Raynaud \cite{Ra}).
\begin{defn} Let $P_{H_n}$ be the {\em generalized Picard scheme} representing
the functor, which associates to $S'\to S$ the isomorphism classes of 
pairs $(\Lh,\alpha)$, where $\Lh$ is a line bundle on $E_n\times_SS'$ and
$\alpha:h_n^*\Lh\isom\Oh_{H_n}$ is a trivialization of $\Lh$ along $H_n\times_SS'$.
\end{defn}
That $H_n$ is a rigidificator in the sense of \cite{Ra}
follows from the fact that $H_n$ contains the section $e_n:S\to E_n$.
Denote by $P_n$ the Picard scheme of $E_n$. Then we have an  exact sequence
\[
0\to T_{H_n}\to P_{H_n}\to P_n\to 0,
\]
where $T_{H_n}$ is the torus with character group $I[H_n]:=\ker(p_{n*}\Z\to \Z)$.
The $l^r$-torsion of $P_{H_n}$ can be identified as follows:
\begin{lemma}\label{p-torsion} There is a canonical isomorphism
\[
R^1\pi_{n!}\Lambda_{r}(1)\isom  P_{H_n}[l^r].
\]
\end{lemma}
\bew Define $_{H_n}\Gm:=\ker(\Gm\to h_{n*}h^*_n\Gm)$. Then
$P_{H_n}=R^1\pibar_*(_{H_n}\Gm)$ (see \cite{SGA4}, Expose XVIII proposition 
1.5.14), 
where the higher direct image is
taken for the flat topology. The sequence
from loc. cit. lemma 1.6.1
\[
0\to j_{n!}\mu_{l^r}\to _{H_n}\Gm\xrightarrow{[l^r]}_{H_n}\Gm\to 0
\]
gives an isomorphism $R^1\pibar_*j_{n!}\mu_{l^r}\isom  P_{H_n}[l^r]$.
\bewende
On $E_n\times_SU_n$ we have the line bundle $\Oh(\Delta_n)$ associated
to the Cartier divisor defined by the diagonal $\Delta_n$. By definition
$\Oh(\Delta_n)$ sits in an exact sequence
\[
0\to \Oh\to \Oh(\Delta_n)\to \Delta_{n*}\Oh_{U_n}\to 0,
\]
which induces a trivialization
of $\Oh(\Delta_n)$ along $H_n\times_SU_n$. Thus we get a section 
$\Delta_n:\Z\to \pi_n^*P_{H_n}$ of \'etale sheaves. Adjunction gives
a map 
\begin{equation}
\Delta_n:p_{n*}\Z\to \pi^*P_{H_n} 
\end{equation}
also denoted by $\Delta_n$ by abuse of notation. 

\subsubsection{Comparison with a one-motive}
Consider
$p_{n*}\Z\to \pi^*P_{H_n}$ as a complex of sheaves in degree $0$ and $1$.
 Note that $p_{n*}\Z\to \pi^*P_{H_n}$ is not a one-motive 
because $P_{H_n}$ is not a semi-abelian scheme. We have an 
exact sequence $0\to  P_{H_n}^0\to P_{H_n}\to \Z\to 0$ and if
we let $I[H_n]$ be the kernel of the composition $p_{n*}\Z\to \Z$ we get
a quasi-isomorphism
\[
[I[H_n]\to  P_{H_n}^0]\isom [p_{n*}\Z\to \pi^*P_{H_n}].
\]
Here $I[H_n]\to  P_{H_n}^0$ is of course a one-motive.

\begin{thm}\label{onemotivecomp}
There is a canonical isomorphism of \'etale sheaves
\[
\underline{H}^0([p_{n*}\Z\to \pi^*P_{H_n} ]\otimes^{\Bbb L}\Z/l^r\Z)\isom
(\polgeom_n)^{\vee}(1)
\]
which is compatible with the morphism 
$(\polgeom_n)^{\vee}(1)\to (\polgeom_m)^{\vee}(1)$ induced by the trace map 
$\Logeom_m\to \Logeom_n$.
\end{thm}

\bew
Let $t:S\to E$ be a section and $H_{t,n}:=p_n^{-1}(t)$ be its
preimage in $E_n$. Let $\rho_n:H_{t,n}\to E_n$ be the
embedding and denote by $g_n$ the open immersion of the complement
\[
g_n:E_n\ohne H_{t,n}\hookrightarrow E_n.
\]
Denote also by $p_{t,n}:H_{t,n}\to S$ the structure map.
We have an exact sequence of \'etale sheaves
\[
0\to \Gm\to g_{n*}g_n^*\Gm\to \rho_{n*}\Z\to 0.
\]
On the other hand let
$_{H_n}\Gm:=\ker(\Gm\to h_{n*} h_n^*\Gm)$,
where $h_n:H_n\to E_n$ is the embedding. Then we get  an exact 
sequence
\[
0\to _{H_n}\Gm\to g_{n*}g_n^*({_{H_n}}\Gm)\to \rho_{n*}\Z\to 0.
\]
The long exact cohomology sequence for $R\pibar_{n*}$ taken for 
the flat topology gives
\[
(p_{t,n})_*\Z\to \pi^*R^1\pibar_{n*} ({_{H_n}}\Gm)\to R^1\pibar_{n*}
g_{n*}g_n^*({_{H_n}}\Gm)\to 0
\]
because $R^1 p_{t,n}\Z=0$. 
We have an identification $R^1\pibar_{n*} ({_{H_n}}\Gm)\isom P_{H_n}$
and it is clear that the map $ (p_{t,n})_*\Z\to R^1\pibar_{n*} ({_{H_n}}\Gm)\isom P_{H_n}$
is injective. We get an isomorphism of complexes
\[
[ (p_{t,n})_*\Z\to  P_{H_n}]\otimes^{\Bbb L}\Z/l^r\Z\isom
[R^1\pibar_{n*}g_{n*}g_n^*({_{H_n}}\Gm)[1]]\otimes^{\Bbb L}\Z/l^r\Z.
\]
Using the sequence 
\[
0\to j_{n!}\mu_{l^r}\to _{H_n}\Gm\xrightarrow{[l^r]}_{H_n}\Gm\to 0
\]
as in the proof of lemma \ref{p-torsion}, we see that
\[
\underline{H}^0([( p_{t,n})_*\Z\to \pi^*P_{H_n} ]\otimes^{\Bbb L}\Z/l^r\Z)\isom
R^1\pibar_{n*}g_{n*}g_n^*(j_{n!}\mu_{l^r}),
\]
where $j_n:U_n\to E_n$.
Let us compute $t^*\polgeom_n$. By base change and the definition \ref{poldefn}
it is a straightforward computation that
\[
t^*\polgeom_n=R^1\pibar_{n*}g_{n!}g_n^*j_{n*}\mu_{l^r}
\]
and by Poincar\'e duality we get 
\[
t^*(\polgeom_n)^{\vee}(1)=R^1\pibar_{n*}g_{n*}g_n^*j_{n!}\mu_{l^r}.
\]
If we apply this to the universal section $\Delta:U\to E\times_SU$ 
we get the desired result.
\bewende
Here is an explicit way to get the extension $(\polgeom_n)^{\vee}(1)$
from $[p_{n*}\Z\to \pi^*P_{H_n}]$. Consider the exact sequence
\[
0\to\pi^* P_{H_n}[l^r]\to \pi^*P_{H_n}\xrightarrow{l^r}
\pi^*P_{H_n}\to \Z/l^r\Z\to 0
\]
defined by the $l^r$-multiplication
and the map $p_{n*}\Z\to \pi^*P_{H_n}$. The pull-back gives an 
extension
\[
0\to\pi^* P_{H_n}[l^r]\to\Eh\to p_{n*}\Z\to \Z/l^r\Z\to 0.
\]
If we tensor this with $\Z/l^r\Z$ over $\Z$, then
we get an exact sequence
\[
0\to\pi^* P_{H_n}[l^r]\to\Eh\otimes_{\Z}\Z/l^r\Z\to 
p_{n*}\Z/l^r\Z\to \Z/l^r\Z\to 0,
\]
 because 
the kernel of $p_{n*}\Z\to \Z/l^r\Z$ is a free $\Z$-module.
With the above theorem we conclude that $(\polgeom_n)^{\vee}(1)\isom\Eh\otimes_{\Z}\Z/l^r\Z$. 

\subsubsection{The class of the geometric polylog}

Recall from section \ref{reformulation} that the polylog
$(\polgeom_{n})^{\vee}(1)$ defines a class in 
\[
\Ext^1_U(I_{n,\Lambda_r}, \pi^*_nR^1\pi_{n!}\Lambda_{r,U}(1)),
\]
where $I_{n,\Lambda_r}$ is the kernel of the map 
$p_{n*}\Lambda_{r,U}\to \Lambda_{r,U}$.
With lemma \ref{p-torsion} we can write 
\[
\cl((\polgeom_{n})^{\vee}(1))\in \Ext^1_U(I_{n,\Lambda_r},P_{H_n}[l^r]).
\]
Let $[\af]$ be an isogeny of $E$ of degree $\Norm\af$ prime to $l$.
We can consider the pull-back $[\af]^*(\polgeom_{n})^{\vee}(1)$ on
$[\af]^{-1}U$. Recall from lemma \ref{translation} that 
we have an isomorphism $p_{n*}\Lambda_{r,U}\isom [\af]^*p_{n*}\Lambda_{r,U}$.
Thus $[\af]^*(\polgeom_{n})^{\vee}(1)$ defines a class in
\[
\Ext^1_{[\af]^{-1}U}(I_{n,\Lambda_r},P_{H_n}[l^r]).
\] 
Similarly consider  the $\Norm\af$-multiplication on 
$(\polgeom_{n})^{\vee}(1)$
restricted to $[\af]^{-1}U$. This gives also  an extension class
in the above group. 
Denote by 
\[
\cl(([\af]^*-\Norm\af)(\polgeom_{n})^{\vee}(1))
\]
 the difference of these two classes. Our aim is to compute the 
specialization of this class. For this we show first, that
this class is in the image of 
\[
\Ext^1_{[\af]^{-1}U}(p_{n*}\Lambda_r,P_{H_n}[l^r])\to
\Ext^1_{[\af]^{-1}U}(I_{n,\Lambda_r},P_{H_n}[l^r]).
\] 
\begin{lemma}  The class of $\cl(([\af]^*-\Norm\af)(\polgeom_{n})^{\vee}(1))$
on $[\af]^{-1}U$ maps to zero under the boundary map
\[
\Ext^1_{[\af]^{-1}U}(I_{n,\Lambda_r}, P_{H_n}[l^r])\to
\Ext^2_{[\af]^{-1}U}(\Lambda, P_{H_n}[l^r]).
\]
\end{lemma}
\bew The sheaf $[\af]^*(\polgeom_{n})^{\vee}(1)$ can be described 
as the pull-back of the sequence
\[
0\to\pi^* P_{H_n}[l^r]\to \pi^*P_{H_n}\xrightarrow{l^r}
\pi^*P_{H_n}\to \Z/l^r\Z\to 0
\]
by the map 
\[
[\af]^*\Delta_n:p_{n*}\Z\to[\af]^*p_{n*}\Z\to [\af]^*\pi^*P_{H_n}\isom \pi^*P_{[\af]^{-1}H_n}\to \pi^*P_{H_n},
\]
 where the last map is forgetting the trivialization
on $[\af]^*H_n\ohne H_n$. The sheaf $\Norm\af(\polgeom_{n})^{\vee}(1)$
is the pull-back by the map $\Norm\af\Delta_n:p_{n*}\Z\to \pi^*P_{H_n}\xrightarrow{\Norm\af}
\pi^*P_{H_n}$ of the same sequence. 
We have to compute the composition  of these maps
with $\pi^*P_{H_n}\to \Z$. This composition factors in both cases 
as $p_{n*}\Z\to \Z\xrightarrow{\Norm\af}\Z$. Thus, the difference is 
zero. This implies the claim.
\bewende
Consider the difference of the two maps in the above proof and
denote this by $([\af]^*-\Norm\af)\Delta_n$. By adjunction this
defines a map 
\[
([\af]^*-\Norm\af)\Delta_n:\Z\to \pi_n^*P_{H_n},
\]
which by the above proof even factors through $ \pi_n^*P_{H_n}^0$.
The line bundle which is associated with this map is easily described:

Recall from the end of section \ref{picardscheme} that the line bundle 
$\Oh(\Delta_n)$ on $E_n\times_SU_n$  with its canonical trivialization
gives a section $\Delta_n:\Z\to \pi_n^*P_{H_n}$.
Then the line bundle 
$[\af]^*\Oh(\Delta_n)\otimes\Oh(\Delta_n)^{\otimes -\Norm\af}$
on $E_n\times_S[\af]^{-1}U_n$
defines the section $([\af]^*-\Norm\af)\Delta_n$.

Denote by $\delta$ the composition
\[
\delta:H^0([\af]^{-1}U_n,T_{H_n})\to H^1([\af]^{-1}U_n,T_{H_n}[l^r])\to H^1([\af]^{-1}U_n,P_{H_n}[l^r]),
\]
where the first map is the boundary map for the short exact sequence
\begin{equation}\label{lrsequence2}
0\to  T_{H_n}[l^r]\to T_{H_n}\xrightarrow{[l^r]}T_{H_n}\to 0.
\end{equation}

\begin{prop}\label{polaslinebundle} The section
$([\af]^*-\Norm\af)\Delta_n$ factors through $T_{H_n}$
and up to sign, there is an equality
\[
\cl(([\af]^*-\Norm\af)(\polgeom_{n})^{\vee}(1))=
\pm \delta(([\af]^*-\Norm\af)\Delta_n)
\]
in $ H^1([\af]^{-1}U_n,P_{H_n}[l^r])$.
\end{prop}
\bew The composition of $([\af]^*-\Norm\af)\Delta_n$ with the 
map $P_{H_n}^0\to P_n^0$ to the Picard scheme is zero, because 
$[\af]^*\Oh(\Delta_n)\isom\Oh(\Delta_n)^{\otimes \Norm\af}$.
This gives the first claim.
The extension class $\delta(([\af]^*-\Norm\af)\Delta_n)$ is given by the 
pull-back of the sequence (\ref{lrsequence2}) with the map 
$([\af]^*-\Norm\af)\Delta_n:\Z\to \pi_n^*T_{H_n}$.
Using the explicit description following the proof of theorem \ref{onemotivecomp} proves the result.
\bewende

\subsection{Computation of the specialization of the polylog}

Let $t:S\to U$ be an $N$-torsion point. We want to compute the 
specialization at $t$ of the polylog.

\subsubsection{Specialization of the class of the geometric polylog}
\label{specializationclass}
Let $t:S\to U$ be an $N$-torsion point and $[\af]:E\to E$ an
isogeny of degree $\Norm\af$ relatively prime to $Nl$. 
Let $H_{n,t}:=p_n^{-1}(t)\subset U_n$ be the inverse image of the
$N$-torsion point $t:S\to U$ in $U_n$.
To compute the specialization of the polylog following section \ref{specsec} we need
to compute the restriction of the class 
$([\af]^*-\Norm\af)(\polgeom_{n})^{\vee}(1)$
in $H^1([\af]^{-1}U_n,T_{H_n}[l^r])$ to $H_{n,t}$. We let $t:H_{n,t}\to [\af]^{-1}U_n$
also be the inclusion of $H_{n,t}$ by abuse of notation.
 Denote the restriction of the polylog by 
\[
\cl(t^*([\af]^*-\Norm\af) (\polgeom_{n})^{\vee}(1))\in 
H^1(H_{n,t},T_{H_n}[l^r] ).
\]
We have to compute the section 
$t^*([\af]^*-\Norm\af)\Delta_n:\Z\to t^*T_{H_n}$, which is defined
by the trivialization of the (trivial) line bundle 
$t^*[\af]^*\Oh(\Delta_n)\otimes t^*\Oh(\Delta_n)^{\otimes -\Norm\af}$
on $H_{n,t}$.

To describe this section in terms of functions, we 
extend our base to a Galois covering $S_n$ of $S$, where 
$H_{n,t}$ is rational and then
use descent: Let $G_n$ be the Galois covering group of
$S_n/S$. Consider $\prod_{h_n\in H_n(S_n)}\Gm$ and write a typical element of 
this product as $\sum_{h_n}g_{h_n}(h_n)$. Similar definitions apply to
$\prod_{h_n\in H_n(S_n)}\mu_{l^r}$.
Let us write $T_{H_n}$ as a quotient:
\begin{lemma}\label{Tasquot} The group  $T_{H_{n}}(S_n)$ can be identified as the 
$G_n$-invariant elements of the quotient
\begin{align*}
\prod_{h_n\in H_{n}(S_{n})}\Gm(S_n)&\to T_{H_{n}}(S_n)\\
\sum_{h_n}(g_{h_n})(h_n)&\mapsto 
\left(\frac{g_{h_n}}{g_{e_n}}\right)(h_n).
\end{align*}
\end{lemma}
\bew
Clear because the character group of $T_{H_n}$ is $\ker(p_{H_n*}\Z\to \Z)$.
\bewende
This gives 
\[
H^0(H_{n,t}\times_SS_n, T_{H_{n}})=\prod_{t_n\in H_{n,t}(S_n)}T_{H_{n}}(S_n).
\]
With this notation we rewrite the section $t^*([\af]^*-\Norm\af)\Delta_n:\Z\to t^*T_{H_n}$.
Note first that $t^*[\af]^*\Oh(\Delta_n)=[\af]^*\Oh([\af]H_{n,t})$ and
$t^*\Oh(\Delta_n)^{\otimes -\Norm\af}=\Oh(H_{n,t})^{\otimes -\Norm\af}$.
The section $t^*([\af]^*-\Norm\af)\Delta_n:\Z\to t^*T_{H_n}$
can now be written as
\[
\sum_{t_n\in H_{n,t}(S_n)}s(t_n),
\]
where $s(t_n)\in T_{H_{n}}(S_n)$. This section is
given by the trivialization of the line bundle 
$[\af]^*\Oh([\af]t_n)\otimes\Oh(t_n)^{\otimes -\Norm\af}$.
 Using the lemma, we write
$s(t_n)=\sum_{h_n}s(t_n)_{h_n}(h_n)$, with $s(t_n)_{h_n}\in \Gm$. 

We want to compute $s(t_n)_{h_n}$ for $h_n\in H_n(S_n)\ohne e_n$.
For this we need the elliptic units on $E_n$.

\subsubsection{Elliptic units}
Elliptic units were introduced and studied by Robert and Gillard.
An algebraic approach to elliptic units was proposed by Kato. 

Let us recall the  definition of the elliptic units
and their characterization in \cite{Ka2} III 1.1.5. (see also \cite{Scho}).

With Kato let us make the definition that  $a,b\in\End_{\Oh_S}(E)$
are relatively prime $(a,b)=1$, if $\ker(a)\cap \ker(b)=e$ where $e$ is
endowed with the reduced subscheme structure and $ab=ba$.
\begin{thm}\label{thetadefn}(see \cite{Ka2} III 1.1.5. and \cite{Scho} thm 1.2.1)  
Let $a\in \End_{\Oh_S}(E)$ be an endomorphism
with $(a,6)=1$. There is a unique 
section 
\[
\theta_a\in\Oh^*(E\ohne \ker a)
\]
compatible with base change in $S$, with the following properties:\\
i)  $\Div(\theta_a)= \deg(a)(e)-\ker a$\\
ii) for any  $b\in\End_{\Oh_S}(E) $ with $(a,b)=1$ 
\[
b_*\theta_a=\theta_a.
\]
iii) Moreover, for  $b\in\End_{\Oh_S}(E) $ with $(6,b)=1$ and $ab=ba$ 
\[
\frac{\theta_a\verk b}{\theta_a^{\deg(b)}}=
\frac{\theta_b\verk a}{\theta_b^{\deg(a)}}.
\]
\end{thm}

\begin{defn} The values of $\theta_a$ at torsion sections 
$t:\Oh_S\to E\ohne \ker a$  are called {\em elliptic units}.
\end{defn}

\subsubsection{The specialization of the polylog in terms of elliptic units}
We turn to the computation of $s(t_n)=\sum_{h_n}s(t_n)_{h_n}(h_n)$
from the end of section \ref{specializationclass}.

\begin{prop} Let $[\af]$ be relatively prime to 
$Nl$. Then the section $s(t_n)_{h_n}\in \Gm(S_n)$ is given by 
\[
s(t_n)_{h_n}=\left(\frac{\theta_{\af}(h_n-t_n)}{\theta_{\af}(-t_n)}\right)^{-1}
\]
for $h_n\in H_n(S_n)\ohne e_n$ using the identification from lemma 
\ref{Tasquot}.
\end{prop}
\bew The section $s(t_n)_{h_n}$ is given by evaluating an isomorphism
\[
[\af]^*\Oh([\af]t_n)\isom\Oh(t_n)^{\otimes \Norm\af}
\]
at $h_n$.  Let $T_{-t_n}$ be the translation
with $-t_n$ on $E_n$, then 
\[
[\af]^*\Oh([\af]t_n)\isom T_{-t_n}^*[\af]^*\Oh(e_n)\mbox{ and }
\Oh(t_n)^{\otimes \Norm\af}\isom T_{-t_n}^*\Oh(e_n)^{\otimes \Norm\af}.
\]
The function $\theta_{\af}$ gives a section of 
$\Oh(e_n)^{\otimes \Norm\af}\otimes [\af]^*\Oh(e_n)^{\otimes -1}$
and thus $T_{-t_n}^*\theta_{\af}$ gives a section of 
\[
\Oh(t_n)^{\otimes \Norm\af}\otimes [\af]^*\Oh([\af]t_n)^{\otimes -1}.
\]
This proves the claim.
\bewende

\begin{rem} Note that in this proposition we do not really need the elliptic
units, because any function with the right divisor would also describe the
sections $s(t_n)$. The point is that there are two more problems to
solve before we get the specialization of the elliptic polylog. The first 
problem is rather trivial. It is the computation of the map $\pr_t$ from
(\ref{prdefn}),
which is nothing but the identification of $H_{n,t}$ with $H_n$ via
the projection map. 
The second problem is that we still have to compute the splitting from
\ref{splitting} and it is here that the elliptic units will be necessary because
of their norm compatibility property.
\end{rem}
\\

Recall the definition of the map $\pr_t$ from \ref{prdefn}. 
We write
\begin{align*}
\pr_t:H_{n,t}(S_n)&\to H_n(S_n)\\
t_n&\mapsto \tilde{t}_n
\end{align*}
for the projection of $t_n$ to $H_n$. Consider the section
\[
\prod_{\tilde{t}_n\in H_n(S_n)}
\sum_{h_n}
\left(\frac{\theta_{\af}(h_n-t_n)}{\theta_{\af}(-t_n)}\right)(h_n)
\]
of $\prod_{\tilde{t}_n\in H_n(S_n)}T_{H_n}(S_n)$. 
It is invariant under the group $G_n$ and thus defines a section
in $H^0(H_n,T_{H_n})$. Denote by 
\[
\delta :H^0(H_n,T_{H_n})\to H^1(H_n,T_{H_n}[l^r])
\]
the boundary map.
The main result of this section can now be formulated as 
follows:
\begin{prop} The class of 
$\pm\pr_tt^*([\af]^*-\Norm\af)(\polgeom_{n})^{\vee}(1)$ is given
by
\[
\delta\left(\prod_{\tilde{t}_n\in H_n(S_n)}
\sum_{h_n}
\left(\frac{\theta_{\af}(h_n-t_n)}{\theta_{\af}(-t_n)}\right)(h_n)\right)\in
H^1(H_n,T_{H_n}[l^r]),
\]
where $t_n$ is an $l^nN$-division point.
\end{prop}

\subsubsection{The splitting}
In the last section we have computed 
\[
\cl(\pr_tt^*([\af]^*-\Norm\af)(\polgeom_{n})^{\vee}(1))\in H^1(H_n,T_{H_n}[l^r]).
\]
Recall  that the torus $T_{H_n}$ is defined by the character group
$I[H_n]:=\ker(p_{H_n*}\Z\to \Z)$. By definition of $\Ihgeom_n$, we get that
$I[H_n]\otimes\Lambda=\Ihgeom_n$ is the augmentation ideal in $\Rhgeom_n$.
In particular
\[
T_{H_n}[l^r]=\underline{\Hom}_S(\Ihgeom_n,\mu_{l^r}).
\]
This gives $H^1(H_n,T_{H_n}[l^r])=H^1(S, \underline{\Hom}_S(\Ihgeom_n,\Rhgeom_n(1)))$.
We have a map  (cf. \ref{map}) 
$\Rhgeom_n(1)\to \underline{\Hom}_S(\Ihgeom_n,\Rhgeom_n(1))$
given explicitly as follows: Write $\Rhgeom_n(1)=\Lambda(1)[H_n]$ and
observe that $\Ihgeom_n\subset \Lambda[H_n]$ is generated by 
$(h_n)-(e_n)$ with $h_n\in H_n$. Then
\begin{align*} 
\Lambda(1)[H_n]&\to \underline{\Hom}_S(\Ihgeom_n,\Lambda(1)[H_n])\\
(h_n)&\mapsto \{ (h'_n)-(e_n)\mapsto ({h'}_n+h_n)-(h_n)\}.
\end{align*}
 This induces a map (cf. corollary \ref{injection})
\begin{equation}\label{splittingmap}
H^1(S,\Rhgeom_n(1))=H^1(H_n,\mu_{l^r})\to H^1(H_n,T_{H_n}[l^r]).
\end{equation}
To compute the splitting it suffices to write down a norm compatible
 element in $H^1(H_n,\mu_{l^r})$, 
which maps to the class of 
$\pm\pr_tt^*([\af]^*-\Norm\af)(\polgeom_{n})^{\vee}(1)$ 
in $H^1(H_n,T_{H_n}[l^r])$.
Recall 
that $\tilde{t}_n$ is the projection to $H_n$. 
\begin{prop}\label{splittingelement} The element
\[
\delta\sum_{\tilde{t}_n\in H_n(S_n)}\theta_{\af}(-t_n)(\tilde{t}_n) \in
H^1(H_n,\mu_{l^r})
\]
maps under (\ref{splittingmap}) to 
$\pm\cl(\pr_tt^*([\af]^*-\Norm\af)(\polgeom_{n})^{\vee}(1))$.
\end{prop}
\bew This is a straightforward computation. The element 
\[
\sum_{\tilde{t}_n}\theta_{\af}(-t_n)(\tilde{t}_n) 
\]
maps to the mapping, which sends $(h_n)-(e_n)$ to 
\[
\sum_{\tilde{t}_n}\theta_{\af}(-t_n)({h}_n+\tilde{t}_n)
-\sum_{\tilde{t}_n}\theta_{\af}(-t_n)(\tilde{t}_n).
\]
The first sum can be rewritten as
\[
\sum_{\tilde{t}_n}\theta_{\af}(h_n-t_n)(\tilde{t}_n)
\]
so that we get $\sum_{\tilde{t}_n}\frac{\theta_{\af}(h_n-t_n)}{\theta_{\af}(-t_n)}(\tilde{t}_n)$.
This is the desired result.
\bewende

We  show now that the element defined in proposition \ref{splittingelement}
is norm compatible if we vary $n$. This will imply that we have actually 
computed the splitting and thus the $l$-adic Eisenstein classes.

Let $\Norm_{n,n'}$  be the norm map from $E_n$ to $E_{n'}$ for $n\ge n'$..

\begin{prop} In $H^1(H_{n'},\mu_{l^r})$ the following equality holds:
\[
\Norm_{n,n'}\delta\sum_{\tilde{t}_n\in H_n(S)}\theta_{\af}(-t_n)
(\tilde{t}_n) =\delta\sum_{\tilde{t}_{n'}\in H_{n'}(S)}
\theta_{\af}(-t_{n'})(\tilde{t}_{n'}) 
\]
where $[l^n]t_n=t$ and $[l^{n-n'}]t_n=t_{n'}$.
\end{prop}
\bew From theorem \ref{thetadefn} we know that 
$\Norm_{n,n'}\theta_{\af}(-t_n)=\theta_{\af}(-t_{n'})$ because
$\af$ is prime to $Nl$. This proves the claim.
\bewende

It is clear that the element 
$\delta\sum_{\tilde{t}_n\in H_n(S_n)}\theta_{\af}(-t_n)(\tilde{t}_n)$
is compatible with the reduction map $\Lambda_r\to\Lambda_{r'}$ for
$r\ge r'$. Hence we get an element
\[
\left(\delta\sum_{\tilde{t}_n\in H_n(S_n)}\theta_{\af}(-t_n)(\tilde{t}_n)\right)_n
\in H^1(S, \Rhgeom_{\Z_l}(1)).
\]
Recall the map 
\[
H^1(S, \Rhgeom_{\Z_l}(1))=\Ext^1_{S,\Z_l}(\Z_l,\Rhgeom_{\Z_l}(1))\xrightarrow{a}
\Ext^1_{S,\Rhgeom_{\Z_l}}(\Ihgeom_{\Z_l},\Rhgeom_{\Z_l}(1))
\]
from corollary \ref{injection}.
\begin{lemma}\label{split} The element $\left(\delta\sum_{\tilde{t}_n\in H_n(S_n)}\theta_{\af}(-t_n)(\tilde{t}_n)\right)_n$
maps to the class of $\polgeom_{\Z_l}$ under the map $a$.
\end{lemma}
\bew This follows immediately from the above proposition.
\bewende

\subsubsection{The main theorem on the specialization of the elliptic
polylog}

It remains to compute the moment map to get the $l$-adic Eisenstein
classes explicitly.

\begin{thm}\label{polastheta} Let $\beta=\sum_{t\in E[N](S)\ohne e}n_t(t)$ and $[\af]:E\to E$
an isogeny relatively prime to $Nl$. Then for $k >0$ the $l$-adic Eisenstein class
\[
(\beta^*([\af]^*-\Norm\af)\pol_{\Q_l})^k\in
H^1(S, \Sym^k\shH_{\Q_l}(1))
\]
is given by 
\[
\pm\frac{1}{k!}\left(\delta\sum_{t\in E[N](S)\ohne e}n_t\sum_{[l^n]t_n=t}
\theta_{\af}(-t_n)\tilde{t}_n^{\otimes k}\right)_n
\]
where $\tilde{t}_n$ is the projection of $t_n$ to $E[l^n]$.
\end{thm}
\bew The recipe to compute the moment map from lemma \ref{momentmap} 
combined with lemma \ref{split} gives immediately the result.
\bewende

\section{Proof of the main theorem}

In this section we will carry out the actual comparison
between the space $r_p(\Rh_{\psi})$ and 
the Soul\'e elements $e_p(\Chbarinfty\otimes T_pE(k))$ defined
by elliptic units.

\subsection{Comparison with the Soul\'e elements}
We first transfer the result from \ref{polastheta} into the setting
of \ref{Beilinson} and then compare these elements with the
Soul\'e map $e_p$.
\subsubsection{The specialization of the elliptic polylog}
The theorem \ref{polastheta} gives us an explicit description
of 
\[
(\beta^*([\af]^*-\Norm\af)\pol_{\Q_l})^{2k+1}H^1_{\et}(\Oh_S, \Sym^{2k+1}\shH_{\Q_p}(1)),
\]
 which we
now translate in the   setting of section \ref{Beilinson}. Note first that $p=l$ and 
$S=\Oh_S$ and $\shH_{\Q_p}=T_pE\otimes\Q_p$. 
Let $\af\in \Oh_K$ be prime to $6p\ff$.
 Let $\theta_{\af}$ be the function defined in \ref{thetadefn}. 
To have shorter 
formulas we introduce the following notation:
Define for $\tilde{t}_r\in E[p^r]$
\[
\gamma(\tilde{t}_{r})^k:=<\tilde{t}_r,
\sqrt{d_K}\tilde{t}_r>^{\otimes k}
\]
where $<\_,\_>$ is the Weil pairing and $\sqrt{d_K}$ is a root of the 
discriminant of $K/\Q$.
Note that $\gamma(\psi(\pf)\tilde{t}_{r})^k=(\Norm \pf)^k\beta(\tilde{t}_{r})^k$.
Recall from the end of section \ref{specsec} that we have an endomorphism 
$([\af]^*-\Norm\af)$ of $H^1_{\et}(\Oh_S, \Sym^{2k+1}\shH_{\Q_p}(1))$,
which is multiplication with $\af^{2k+1}\Norm\af-\Norm\af$. 
\begin{thm}\label{rpdescription}
Let $p\nmid 6\Norm\ff$ and denote for a $p^r\Norm\ff$-torsion point
$t_r$ by $\tilde{t}_r$ its projection to $E[p^r]$. Then with $t=\Omega f^{-1}$
\[
([\af]^*-\Norm\af)r_p(\xi)=\pm\frac{\Norm\ff^{3k+2}L_p(\bar{\psi},-k)^{-1}}{2^{k-1}\psi(f)}
\left(\delta\Norm_{K(\ff)/K} \sum_{p^rt_r=t}\theta_{\af}(-t_{r})\otimes \tilde{t}_r
\otimes \gamma(\tilde{t}_r)\right)_r
\]
where $\tilde{t}_r$ is the projection of $t_r$ to $E[p^r]$.
\end{thm}
\bew We have by definition and 
 and theorem \ref{eiseq} 
\begin{align*}
r_p(\xi)&=\frac{(-1)^{k-1}(2k+1)!L_p(\bar{\psi},-k)^{-1}}
{2^{k-1}\psi(f)\Norm \ff^{k}}\Kh_{\Mh}r_p(\Eh^{2k+1}_{\Mh}(\beta))\\
&=\frac{(-1)^{k-1}(2k+1)!\Norm \ff^{3k+2}L_p(\bar{\psi},-k)^{-1}}
{2^{k-1}\psi(f)}\Kh_{\Mh}(\beta^*\pol_{\Q_l})^{2k+1}.
\end{align*}
By definition of $[\af]^*\pol_{\Q_l}$ we see immediately, that
\[
([\af]^*-\Norm\af)(\beta^*\pol_{\Q_l})^{2k+1}=(\beta^*([\af]^*-\Norm\af)\pol_{\Q_l})^{2k+1}.
\]
With the above notation, we have
\[
\Kh_{\Mh}( \tilde{t}_r^{\otimes 2k+1})= \tilde{t}_r\otimes\gamma(\tilde{t}_r)^k
\]
so that theorem \ref{polastheta} gives:
\[
([\af]^*-\Norm\af)r_p(\xi)=
\pm\frac{\Norm\ff^{3k+2}L_p(\bar{\psi},-k)^{-1}}{2^{k-1}\psi(f)}
\left(\delta\Norm_{K(\ff)/K} \sum_{p^rt_r=t}\theta_{\af}(-t_{r})\otimes \tilde{t}_r
\otimes \gamma(\tilde{t}_r)\right)_r.
\]
This is the desired result.
\bewende

\subsubsection{The comparison theorem}
We want to rewrite the formula in theorem \ref{rpdescription}
in terms of the norm map for $K_n(\ff)/K$.
 
Fix a prime $\pf$ of $K$ where $E$ has good reduction. Define 
a uniformizer by $\pi:=\psi(\pf)$. Denote by
\[
H_{r,t}^{\pf}:=\left\{ t_r\in E[\pf^r\ff]|\pi^rt_r=t\right\}.
\]
We write $t_r=(\tilde{t}_r,\pi^{-r}t)\in E[\pf^r\ff]=E[\pf^r]\oplus E[\ff]$. 
 Denote by $K(\pf^r\ff) $ the ray class
field for $\pf^r\ff$. This is the field where the $ E[\pf^r\ff]$-points are
rational. Let $\sigma_{\pf}$ be the Frobenius at $\pf$ in the 
Galois group of $K(\ff)/K$, then $t_r=(\tilde{t}_r, t^{\sigma_{\pf}^{-r}})$.
Recall that $\gamma(\tilde{t}_r)^k:=<\tilde{t}_r, \sqrt{d_K}\tilde{t}_r>^{\otimes k}$.
Define a filtration of $H_{r,t}^{\pf}$ as follows: 
\[
F^i_{r,t}:=\left\{  t_r=(\tilde{t_r},\pi^{-r}t)\in H_{r,t}^{\pf}
|\pi^{r-i}\tilde{t_r}=0
\right\}.
\]
Thus
\[
H_{r,t}^{\pf}=F^0_{r,t}\supset\ldots\supset F^r_{r,t}=0.
\]
Define $T_{\pf}E:=\prolim_nE[\pf^n]$. 
\begin{thm}\label{thetarel}
Let $\pf$ be as above and
$t_r=(\tilde{t_r},\pi^{-r}t)\in F^0_{r,t}\ohne F^{1}_{r,t}$.
Let $L_{\pf}(\bar{\psi},-k)$ be
the Euler factor for $\bar{\psi}$ at $\pf$ evaluated at $-k$, then
\begin{align*}
L_{\pf}(\bar{\psi},-k)^{-1}&\left(\Norm_{K(\ff)/K}\sum_{s_{r}\in H_{r,t}^{\pf} }\theta_{\af}(-s_{r})
\otimes\tilde{s}_{r}\otimes \gamma(\tilde{s}_{r})^k\right)_r=\\
&\left(\Norm_{K(\pf^r\ff)/K}(\theta_{\af}(-t_{r})
\otimes \tilde{t}_{r}\otimes \gamma(\tilde{t}_{r})^k)\right)_r
\end{align*}
in $H^1(\Oh_S, T_{\pf}E(k+1)\otimes\Q_p)$
for all $\af$  relatively prime to $\pf\ff$.
\end{thm}

\bew Observe that we identified $\Hom_{\Oh_p}(T_pE,\Oh_p)\isom T_pE(-1)$
where  $T_pE$ has now the conjugate linear $\Oh_p$-action.
In particular, $\bar{\psi(\pf)}t_r= t_{r-1}$ for $t_r\in E[\pf^r]$. 
We compute
\begin{align*}
\left(\frac{\bar{\psi(\pf)}}{\Norm \pf^{-k}}\right)^i\Norm_{K(\pf^r\ff)/K(\pf^{r-i}\ff)}
&\left(\theta_{\af}(-t_r)\otimes \tilde{t}_r\otimes 
\gamma(\tilde{t}_{r})^k\right)
=\\
&=\Norm_{K(\pf^r\ff)/K(\pf^{r-i}\ff)}\left(
\theta_{\af}(-(\tilde{t}_r,\pi^{-r}t))\otimes \bar{\psi(\pf)}\tilde{t}_r\otimes 
\gamma(\bar{\psi(\pf)}\tilde{t}_{r})^k\right)\\
&= \left(\Norm_{K(\pf^r\ff)/K(\pf^{r-i}\ff)}\theta_{\af}(-(\tilde{t}_r,\pi^{-r}t))\right)
\otimes\tilde{t}_{r-i}\otimes \gamma(\tilde{t}_{r-i})^k\\
&= \theta_{\af}(-(\tilde{t}_{r-i},\pi^{i-r}t)))\otimes\tilde{t}_{r-i}\otimes \gamma(\tilde{t}_{r-i})^k,
\end{align*}
where we used the distribution relation for $\theta_{\af}$ 
(see \cite{deSh} II 2.5)
\[
\Norm_{K(\pf^r\ff)/K(\pf^{r-i}\ff)}\theta_{\af}(-{t}_r)
=\theta_{\af}(\pi^i{t}_r))
\]
for the last equality.

As the Galois group of $K(\pf^{r-i}\ff)/K(\ff)$ acts simply transitively 
on $ F^{i}_{r,t}\ohne  F^{i+1}_{r,t} $, we get
\begin{align*}
\left(\frac{\bar{\psi(\pf)}}{\Norm \pf^{-k}}\right) 
^i\Norm_{K(\pf^r\ff)/K(\ff)}&
\left(\theta_{\af}(-t_r)\otimes \tilde{t}_r\otimes 
\gamma(\tilde{t}_{r})^k\right)
=\\
&\sum_{t_{r-i}\in F^{i}_{r,t}\ohne  F^{i+1}_{r,t}} \theta_{\af}(-(\tilde{t}_{r-i},\pi^{i-r}t))\otimes\tilde{t}_{r-i}\otimes \gamma(\tilde{t}_{r-i})^k.
\end{align*}
We have $ \theta_{\af}(-(\tilde{t}_{r-i},\pi^{i-r}t)))=
 \theta_{\af}(-(\tilde{t}_{r-i},\pi^{-r}t)))^{\sigma_{\pf}^{i}}$ 
which gives 
\begin{align*}
\left(\frac{\bar{\psi(\pf)}}{\Norm \pf^{-k}}\right)^i\Norm_{K(\pf^r\ff)/K}&
\left(\theta_{\af}(-t_r)\otimes \tilde{t}_r\otimes 
\gamma(\tilde{t}_{r})^k\right)
=\\
&\Norm_{K(\ff)/K}\sum_{t_{r-i}\in F^{i}_{r,t}\ohne  F^{i+1}_{r,t}} 
\left(\theta_{\af}(-(\tilde{t}_{r-i},\pi^{-r}t))\otimes\tilde{t}_{r-i}\otimes \gamma(\tilde{t}_{r-i})^k\right)
\end{align*}
because the norm $\Norm_{K(\ff)/K}$ is the sum over all the
Galois translates, which act trivially on $\tilde{t}_{r-i}$.
If we finally take the sum over $i$ and let $r$ get bigger and bigger
we get
\begin{align*}
L_{\pf}(\bar{\psi(\pf)},-k)&\left(\Norm_{K(\pf^r\ff)/K}
\theta_{\af}(-t_r)\otimes \tilde{t}_r\otimes 
\gamma(\tilde{t}_{r})^k\right)_r=\\
&\left(\Norm_{K(\ff)/K}\sum_{t_{r}\in H_{r,t}^{\pf}}
 \theta_{\af}(-{t}_{r})\otimes\tilde{t}_{r}\otimes 
\gamma(\tilde{t}_{r})^k\right)_r,
\end{align*}
where we used
\[
\sum_{i\ge 0}\left(\frac{\bar{\psi(\pf)}}{\Norm \pf^{-k}}\right)^i =\frac{1}{1-
\frac{\bar{\psi(\pf)}}{\Norm \pf^{-k}}}.
\]
This is the desired result.
\bewende

With theorem  \ref{rpdescription} we get:
\begin{cor}\label{rpasnorm} With the notations of theorem \ref{rpdescription}
\[
([\af]^*-\Norm\af) r_p(\xi)=\pm\frac{\Norm\ff^{3k+2}}{2^{k-1}\psi(f)}\delta
\left(\Norm_{K(p^r\ff)/K} \theta_{\af}(-t_{r})\otimes \tilde{t}_r
\otimes \gamma(\tilde{t}_r)^k\right)_r,
\]
where $p^rt_r=t$ and $t_r$ is a primitive $p^r\ff$-torsion point.
\end{cor}
\bew If $p$ is inert of prime this is just a reformulation of 
theorem \ref{thetarel}. If $p$ is split, $r_p$ decomposes into a
direct sum for the $\pf$ and the $\pf^*$ part. 
Putting them together gives the result. 
\bewende

\subsection{End of proof of the main theorem}
Here we finich the proof of theorem \ref{mainthm} by 
computing the image of the Soul\'e map $e_p$. In the last section we
prove that $r_p$ is injective on $\Rh_{\psi}$ if $H^2_p$ is finite.
\subsubsection{Relation to elliptic units}
Our aim is to show that the elements
\[
\left(\Norm_{K(p^r\ff)/K}\theta_{\af}(t_{r})
\otimes \tilde{t}_{r}\otimes \gamma(\tilde{t}_{r})^k\right)_r
\]
where $\af$ is prime to $6p\ff$
generate $(\Chbarinfty^{\chi}\otimes T_pE(k))_{\Gamma}$,
where $\chi$ is the representation of $\Delta$ on $\Hom_{\Oh_p}(T_pE(k),\Oh_p)$.
\begin{prop} Let $p\nmid 6\Norm\ff$ and
$\af$ be an ideal in $\Oh_p$, which is prime to
$6p\ff$ and such that $\Norm\af\not\equiv 1(\mbox{mod }p)$.
Then the $\Oh_p[[\Gamma]]$-module 
\[
\Chbarinfty^{\chi}\otimes_{\Oh_p} T_pE(k)
\]
is generated by $\left(\theta_{\af}(t_{r})
\otimes \tilde{t}_{r}\otimes \gamma(\tilde{t}_{r})^k\right)_r$, where
$t_r$ is a primitive $p^r\ff$-division point.
\end{prop}
\bew Let $\bff $ be another ideal 
prime to $6p\ff$. Then by theorem \ref{thetadefn} 
\begin{align*}
(\sigma_{\af}-\psi(\af)(\Norm\af)^{k+1})(\theta_{\bff}(t_{r})
\otimes \tilde{t}_{r}\otimes \gamma(\tilde{t}_{r})^k)
&=
\psi(\af)(\Norm\af)^{k}(\theta_{\bff}(t_{r})^{\sigma_{\af}-\Norm\af}\otimes \tilde{t}_{r}\otimes \gamma(\tilde{t}_{r})^k)\\
&=\psi(\af)(\Norm\af)^{k}(\theta_{\af}(t_{r})^{\sigma_{\bff}-\Norm\bff}
\otimes \tilde{t}_{r}\otimes \gamma(\tilde{t}_{r})^k).
\end{align*}
It is enough to show that $\sigma_{\af}-\psi(\af)(\Norm\af)^{k+1}$
is invertible in $\Lambda=\Oh_p[[\Gamma]]$, because $\Chbarinfty^{\chi}$
is a torsion free $\Lambda$-module. $\Lambda$ is a local ring if 
$p$ is inert or prime in $K$ and a product of local rings if $p$ is 
split.  We have $\Lambda/\mf=E[p](k)$, where $\mf$ is either
the maximal ideal or the product of the maximal ideals.
The element $\sigma_{\af}$
acts via $\psi(\af)(\Norm\af)^k$ on $E[p](k)$ and thus 
$\sigma_{\af}-\psi(\af)(\Norm\af)^{k+1}$ is invertible in $\Lambda$
if $\Norm\af\not\equiv 1(\mbox{mod }p)$.
It remains to see that $\gamma(\tilde{t}_{r})^k$ generates
$\Z_p(k)$. We have 
\begin{align*}
<\tilde{t}_r, \sqrt{d_K}\tilde{t}_r>^{\pm 1}&=
\exp(p^{-r}|\Omega|^2(\bar{\sqrt{d_K}}-\sqrt{d_K}))\\
 &= \exp(-2 ip^{-r}|\Omega|^2{\sqrt{|d_K|}})\\
&=\exp(-4\pi ip^{-r}),
\end{align*}
which is for $p\neq 2$ a primitive root of unity.
\bewende

\begin{cor} \label{rpequ}The image $r_p(\Rh_{\psi})$ in $H^1(\Oh_S, T_pE(k+1)\otimes\Q_p)$ 
coincides with the image of 
\[
e_p((\Chbarinfty^{\chi}\otimes T_pE(k))_{\Gamma}).
\]
\end{cor}
\bew As $\frac{\Norm\ff^{3k+2}}{2^{k-1}\psi(f)}$ is prime 
to $p$, this follows from corollary \ref{rpasnorm} and the definition of 
$e_p$ in \ref{souleelements}.
\bewende

To conclude the proof of theorem \ref{mainthm} it remains to 
see the following lemma:

\begin{lemma}The canonical map 
\[
(\Chbarinfty\otimes  T_pE(k))\otimes^{\Lbb}_{\Oh_p[[\Gh]]}\Oh_p\to
( \Chbarinfty\otimes  T_pE(k))_{\Gh}\isom( \Chbarinfty^{\chi}\otimes  T_pE(k))_{\Gamma} 
\]
is an isomorphism and $( \Chbarinfty^{\chi}\otimes  T_pE(k))_{\Gamma}\isom\Oh_p$.
\end{lemma}
\bew By \cite{Ru1} theorem 7.7 we have an isomorphism
\[
\Chbarinfty^{\chi}\isom \Lambda^{\chi}=\Oh_p[[\Gamma]].
\]
This implies that the $\Oh_p[[\Gamma]]$-module is induced
and hence as $\Oh_p[[\Gamma]]$-module isomorphic to
$\Oh_p[[\Gamma]]$. This implies 
$( \Chbarinfty^{\chi}\otimes  T_pE(k))_{\Gamma}\isom\Oh_p$
and the claim of the corollary, because the higher Tor-terms vanish.
\bewende
We get as a corollary part b) of theorem \ref{mainthm}:
\begin{cor} The map 
\[
\Rh_{\psi}\otimes\Z_p\to R\Gamma(\Oh_S,T_pE(k+1)\otimes\Q_p)[1]
\]
induced by $r_p$, gives an isomorphism
\[
{\det}_{\Oh_p}\Rh_{\psi}\isom {\det}_{\Oh_p}R\Gamma(\Oh_S,T_pE(k+1))^{-1}.
\]
\end{cor}
\bew The complex $\Rh_{\psi}\otimes\Z_p\to R\Gamma(\Oh_S,T_pE(k+1)\otimes\Q_p)[1]$ is
isomorphic to 
\[
( \Chbarinfty^{\chi}\otimes  T_pE(k))_{\Gamma}
\xrightarrow{e_p}  R\Gamma(\Oh_S,T_pE(k+1)\otimes\Q_p)[1]
\]
because by \ref{rpequ} $r_p$ and $e_p$ have the same image and 
as $\Oh_p$-modules $( \Chbarinfty^{\chi}\otimes  T_pE(k))_{\Gamma}\isom\Oh_p$
and $\Rh_{\psi}\otimes\Z_p\isom \Oh_p$. Theorem \ref{cohdescr} implies then
the claim.
\bewende

\subsubsection{Finiteness of $H^2_p$ and injectivity of $r_p$}\label{h2finite}
Here we prove the addition of theorem \ref{mainthm}, that the
finiteness of $H^2_p$ implies that $r_p$ is injective on $\Rh_{\psi}$.

Recall that by corollary \ref{rpequ} the image of $r_p$ coincides
with the image of $e_p$. It suffices for the injectivity of
$r_p$ to prove that $r_p(\Rh_{\psi})$ is non zero, because 
$\Rh_{\psi}\isom \Oh_K$.
\begin{prop}
Let $H^2_p$ be finite, then $e_p$ is injective.
\end{prop}

\bew We chow first that $H^2_p$ finite implies the finiteness of
$(\Ahinfty\otimes T_pE(k))_{\Gh}$. Note that this is the cokernel of
\[
(\Uhinfty\otimes T_pE(k))_{\Gh}\to (\Xhinfty\otimes T_pE(k))_{\Gh}.
\]
Computing up to finite groups we get from corollary \ref{complexes}
(using corollary \ref{uinftyisom}) that this cokernel is isomorphic
(up to finite groups) to the cokernel of 
\[
H^1(K\otimes\Q_p, E[p^{\infty}](-k))^*\to
H^1(\Oh_{S_p}, E[p^{\infty}](-k))^*
\]
which is contained in $H^2(\Oh_{S_p},T_pE(k+1))$. This group is 
of course finite if $H^2_p=H^2(\Oh_{S},T_pE(k+1))$ is finite.
Thus $(\Ahinfty\otimes T_pE(k))_{\Gh}$ is finite.
Using lemma 6.2. from \cite{Ru1} we see that this implies that
$(\Ahinfty\otimes T_pE(k))^{\Gh}$ is finite. We will now show that this 
last group controls the kernel of $e_p$. It suffices to show that the
kernel of $e_p$ on $(\Ehbarinfty\otimes T_pE(k))_{\Gh}$ is finite because
by \cite{Ru1} 7.8. both $\Ehbarinfty$ and $\Chbarinfty$ are $\Lambda$-modules
of rank $1$ with $\Ehbarinfty/\Chbarinfty$ a torsion module. 
So suppose that the image of $(\Ehbarinfty\otimes T_pE(k))_{\Gh}$
under $e_p$ has not rank $1$, i.e. is finite. Then, because
$(\Uhinfty\otimes T_pE(k))_{\Gh}\isom H^1(K\otimes \Q_p,E[p^{\infty}](-k))^* $
the image of $(\Ehbarinfty\otimes T_pE(k))_{\Gh}$ in 
$(\Uhinfty\otimes T_pE(k))_{\Gh}$ must be finite as well. 
The kernel of the map 
\[
(\Ehbarinfty\otimes T_pE(k))_{\Gh}\to (\Uhinfty\otimes T_pE(k))_{\Gh}
\]
is $H_1(\Gh,\Uhinfty/\Ehbarinfty\otimes T_pE(k))$ (group homology). 
On the other hand, up to finite groups, corollary \ref{complexes}
implies that $H_1(\Gh, \Xhinfty\otimes T_pE(k))\isom H^2(\Oh_{S_p}, E[p^{\infty}](-k))^*$.
By lemma \ref{exactseq} we get a commutative diagram (up to finite groups)
\[
\begin{CD}
H_1(\Gh,\Uhinfty/\Ehbarinfty\otimes T_pE(k))@>\alpha>>
H_1(\Gh, \Xhinfty\otimes T_pE(k)\\
@VVV@VVV\\
(\Ehbarinfty\otimes T_pE(k))_{\Gh}@>e_p>>H^1(\Oh_{S_p}, T_pE(k+1))\\
@VVV@VVV\\
(\Uhinfty\otimes T_pE(k))_{\Gh}@>\isom>>H^1(K\otimes \Q_p,E[p^{\infty}](-k))^*.
\end{CD}
\]
The kernel of the map $\alpha $ is a quotient of $(\Ahinfty\otimes T_pE(k))^{\Gh}$
which by the above is finite. Thus, we arrive at a contradiction and
$e_p$ can not be zero on the free part of $(\Ehbarinfty\otimes T_pE(k))_{\Gh}$.
Hence, $e_p$ is non zero on $(\Chbarinfty\otimes T_pE(k))_{\Gh}$.
\bewende

\noindent Guido Kings\\
Mathematisches Institut\\
Westf\"alische Wilhelms-Universit\"at M\"unster\\
Einsteinstr. 62\\
48149 M\"unster\\
Germany\\
e--mail: kings@math.uni-muenster.de


\begin{thebibliography}{99999999}
\bibitem[Ar-Ma]{Ar-Ma}M. Artin, B. Mazur: Etale Homotopy, Lecture 
Notes in Mathematics 100, Springer 1986
\bibitem[Be1]{Be}A. Beilinson: Higher regulators and values of 
L--functions. J. Sov. Math. {\bf 30}, 2036--2070 (1985) 
\bibitem[Be2]{Be1}A. Beilinson: Higher regulators of modular 
curves, Contemporary Math. Vol. 55\,I (1986)
\bibitem[Be3]{Be2} A. Beilinson: Polylogarithm and cyclotomic elements,
manuscript, no date.
\bibitem[Be-Le]{Be-Le}A.\ Beilinson, A.\ Levin: The elliptic polylogarithm, in: U.\ Jannsen et al.(eds.): Motives, Proceedings Seattle 1991,
Providence, RI: American Mathematical Society,  Proc. Symp. Pure Math. 55, Pt. 2, 123-190 (1994)
\bibitem[Bour]{Bour}N. Bourbaki: Groupes et Algebres des Lie,
Hermann, (1972)
\bibitem[Bl-Ka]{Bl-Ka}S. Bloch, K. Kato: L-functions and Tamagawa numbers
of motives, in: P. Cartier et al. eds.: The Grothendieck Festschrift Vol. I,
Birkh\"auser (1990)
\bibitem[Del1]{Del2} P. Deligne: Th\'eorie de Hodge III, Publ. Math. 
IHES, 5--77 (1974)
\bibitem[Del2]{Del1} P. Deligne:  Le groupe fondamental de la droite
projective moins trois points. in: Ihara et al. (eds.): Galois
groups over $\Q$, MSRI Publication (1989)
\bibitem[Den1]{Den} C. Deninger: Higher regulators and Hecke L--series
of imaginary quadratic fields I. Invent. math. {\bf 96} ,1--69 (1989)
\bibitem[Den2]{Den1} C. Deninger: Extensions of motives associated to 
symmetric powers of elliptic curves and to Hecke characters of 
imaginary quadratic fields, in: F.Catanese (ed.):
Arithmetic Geometry, Cortona 1994.
\bibitem[deSh]{deSh}E. deShalit: Iwasawa Theory of Elliptic Curves with
Complex Multiplication, Perspectives in Mathematics vol. 3, Academic Press
(1987)
\bibitem[Hu-Ki1]{Hu-Ki1}A.\ Huber, G.\ Kings: Dirichlet motives via modular curves, Ann. Sci. ENS, {\bf 32}, 313--345 (1999)
\bibitem[Hu-Ki2]{Hu-Ki}A.\ Huber, G.\ Kings: Degeneration of $l$-adic Eisenstein classes and of the elliptic poylog, Invent. math. {\bf 135}, 545--594 (1999)
\bibitem[HuW]{HW}A.\ Huber, J. Wildeshaus: Classical motivic polylogarithm
according to Beilinson and Deligne, Doc. Math. J. DMV {\bf 3} , 27--133 (1998)
\bibitem[Ja1]{Ja1}U. Jannsen: Continous \'etale cohomology, Math. Ann. 
{\bf 280}, 207-245 (1988)
\bibitem[Ja2]{Ja2}U. Jannsen: On the $l$-adic cohomology of varieties
over number fields and its Galois cohomology, in: Ihara et al. (eds.): Galois
groups over $\Q$, MSRI Publication (1989)
\bibitem[Ka1]{Ka1}K. Kato: Lectures on the approach to Iwasawa theory for
Hasse-Weil L-functions via $B_{dR}$, in: J.-L. Colliot-Th\'el\`ene et al.:
Arithmetic Algebraic Geometry, LNM 1553, Springer (1993)
\bibitem[Ka2]{Ka2}K. Kato: Iwasawa theory and p-adic Hodge theory, 
Kodai math. J. {\bf 16} 1--31 (1993)
\bibitem[Ki]{Ki}G. Kings: K-theory elements for the polylogarithm of abelian
schemes, J. reine angew. Math. {\bf 517}, 103--116 (1999) 
\bibitem[Kn-Mu]{Kn-Mu}F. Knudsen, D. Mumford: The projectivity of the 
moduli space of stable curves I: Preliminaries on ``det'' and ``Div'',
Math. Scand. {\bf 39}, 19--55 (1976)
\bibitem[Ra]{Ra}M. Raynaud: Sp\`ecialisation du Foncteur de Picard, 
Publ. Math. IHES, {\bf 38} 27--76 (1970)
\bibitem[Ru1]{Ru2}K. Rubin: Tate Shafarevich groups and L-functions of 
elliptic curves with complex multiplication, Invent. math. {\bf 89}, 527--560
(1987)
\bibitem[Ru2]{Ru1}K. Rubin: The ``main conjectures'' of Iwasawa theory for
imaginary quadratic fields, Invent. math. {\bf 103}, 25--68 (1991)
\bibitem[Ru3]{Ru3}K. Rubin: Elliptic curves with complex multiplication and
the conjecture of Birch and Swinnerton-Dyer, in: Arithmetic Theory of 
Elliptic Curves, Lecture Notes in Mathematics 1716, Springer 1999
\bibitem[Scho]{Scho} A.J. Scholl: An introduction to Kato's Euler system,
in: A.J. Scholl, R.L. Taylor (eds.): Galois representations in 
Arithmetic Algebraic Geometry, Cambridge University Press (1998)
\bibitem[Sch1]{Sch1}P. Schneider: \"Uber gewisse Galoiscohomologiegruppen,
Math. Z. {\bf 168} (1979)
\bibitem[Sch2]{Sch}P. Schneider: Introduction to the
Beilinson conjectures. in: M. Rapoport et al.: Beilinson's
conjectures on special values of $L$--functions. Academic
Press (1988)
\bibitem[Se]{Se}J.-P. Serre: Groupes alg\'ebriques et corps des classes,
Hermann (1959)
\bibitem[Si]{Si}J.H. Silverman: Advanced Topics in the Arithmetic of 
Elliptic Curves, Graduate Texts in Math. vol. 151, Springer (1994)
\bibitem[So1]{So} C. Soul\'e: The rank of \'etale cohomology of varieties
over $p$-adic or number fields, Comp. Math. {\bf 53}, 113--131 (1984)
\bibitem[So2]{So2}C. Soul\'e: p-adic K-theory of elliptic curves,
Duke math. Jour. {\bf 54}, 249--269 (1987)
\bibitem[Wi]{Wi}J.\ Wildeshaus: Realizations of Polylogarithms, Lecture Notes
in Mathematics 1650, Springer 1997.
\bibitem[Win]{Win}K. Wingberg: On the \'etale K-theory of an elliptic curve
with complex multiplication for regular primes, Canad. Math. Bull. {\bf 33},
145--150 (1990)
\bibitem[EGA II]{EGA}A. Grothendieck: \'El\'ements de G\'eom\'etrie Alg\'ebrique II, 
Publ. Math. IHES, {\bf 8} (1961)
\bibitem[SGA4,III]{SGA4} S\'eminaire de G\'eom\'etrie Alg\'ebrique 4,
Th\'eorie des topos et cohomologie \'etale des sch\'emas, Springer LNM 305
(1972)
\bibitem[SGA41/2]{SGA41/2} S\'eminaire de G\'eom\'etrie Alg\'ebrique 4$\frac{1}{2}$, Cohomologie \'etale, Springer LNM 569 (1977) 
\end{thebibliography}
\end{document}